\documentclass{amsart}

\usepackage{amsfonts, amsmath,amscd, amssymb, latexsym, mathrsfs, stmaryrd, verbatim, wasysym }
\usepackage{graphicx}
\usepackage[all]{xy}
\newtheorem{theorem}{Theorem}
\newtheorem{definition}{Definition}

\newtheorem{lemma}[theorem]{Lemma}
\newtheorem{proposition}[theorem]{Proposition}
\newtheorem{remark}[theorem]{Remark}

\numberwithin{equation}{section}
\numberwithin{theorem}{section}
\numberwithin{definition}{section}

 \usepackage{color}
\definecolor{darkgreen}{cmyk}{1,0,1,.2}
\definecolor{m}{rgb}{1,0.1,1}
\DeclareMathOperator{\End}{End}    

\DeclareMathOperator{\Ind}{Ind}

\newcommand{\forget}[1]{}

\def  \nuint {\raise10pt\hbox{$\nu$}\kern-6pt\int}

\newcommand\we{\widehat{E}}

\newcommand\Tr{\operatorname{Tr}}

\def \A{\mathcal A}

\def \F{\mathcal F}
\def \I{\mathcal I}

\newcommand\E{\mathcal E}

\def \Sp {{\cal S}}
\newcommand\B{\mathcal B}

\def \J{\mathcal J}
\def \H {{\cal H}}

\def\Id{{\rm Id}}

\newcommand\cyl{\operatorname{cyl}}

\renewcommand\Im{\operatorname{Im}}

\newcommand\D{\mathcal D}
\newcommand\Di{D\kern-6pt/}
\newcommand\cDi{{\mathcal D}\kern-6pt/}
\newcommand\spi{S\kern-6pt/}
\newcommand \cspi{\Sp\kern-6pt/}

\newcommand\CC{\mathbb C}

\def \cal {\mathcal}

\newcommand\KK{\mathbb K}

\newcommand\RR{\mathbb R}
\newcommand\ZZ{\mathbb Z}

\newcommand\pa{\partial}
\newcommand\Ker{\operatorname{Ker}}

\def\tV{{\tilde V}}
\def\tN{{\tilde N}}

\def\tM{{\tilde M}}
\def\tm{{\tilde m}}

{\catcode`@=11\global\let\c@equation=\c@theorem}

\allowdisplaybreaks[2]
\date{}
 \usepackage{color}
\definecolor{darkgreen}{cmyk}{1,0,1,.2}
\definecolor{m}{rgb}{1,0.1,1}

\title[Eta cocycles]{Eta cocycles}

\author{Hitoshi Moriyoshi}
\address{Graduate School of Mathematics, Nagoya University}
\email{moriyosi@math.nagoya-u.ac.jp}
\author{Paolo Piazza}
\address{Dipartimento di Matematica, Sapienza Universit\`a di Roma} 
\email{piazza@mat.uniroma1.it}



 \dedicatory{\large{Dedicated to  Henri Moscovici on the occasion
  of his sixty-fifth anniversary}}

\begin{document}
\maketitle

\begin{abstract}
We describe a Godbillon-Vey  index formula for longitudinal Dirac operators on a foliated bundle $(X,\F)$ with boundary;
in particular, we define a {\it Godbillon-Vey eta invariant} on $(\pa X,\F_{\pa})$, that is, a secondary
invariant for longitudinal Dirac operators on type III foliations.
Our theorem generalizes the classic Atiyah-Patodi-Singer index formula for $(X,\F)$.
Moreover, employing the Godbillon-Vey index as a pivotal example,
we explain a new approach to higher index theory on geometric
structures with boundary. This is heavily based on
the interplay between the absolute and relative pairing of $K$-theory and
cyclic cohomology for an exact sequence of Banach algebras,
 which in the present
context takes the form
$0\to \mathbf{ \mathfrak{J}} \to \mathbf{ \mathfrak{A}} \to
\mathbf{ \mathfrak{B}} \to 0$
with $ \mathbf{ \mathfrak{J}}$  dense
and holomorphically closed in $C^* (X,\F)$ and   $ \mathbf{ \mathfrak{B}}$
depending only on boundary data.
\end{abstract}

\section{Introduction}
Connes' index theorem for $G$-proper manifolds \cite{Co}, with $G$
an \'etale groupoid, unifies under a single statement most of the
existing (longitudinal) index theorems. We shall focus on a
particular case of such a theorem, that of foliated bundles.
Thus, let $N$ be a closed compact manifold. Let
$\Gamma\to \tN \to N$ be  a Galois $\Gamma$-cover. Let
$T$ be a smooth oriented compact manifold with an action of $\Gamma$
which is assumed to be by diffeomorphisms,
orientation preserving and locally faithful, as in \cite{MN}.
Let $Y=\tN\times_\Gamma T$ and let  $(Y,\mathcal{F})$ be the associated foliation.
(This is an example of $G$-proper
manifold with $G$ equal to the groupoid $T\rtimes \Gamma$.) Let $D$ be a $\Gamma$-equivariant
family of Dirac operators on the fibration $\tN\times T\to T$; such a family induces a longitudinal
Dirac operator on $(Y,\mathcal{F})$.

 If
$T=$ point and $\Gamma=\{1\}$ we have a compact manifold and
Connes' index theorem reduces to the Atiyah-Singer index theorem.
 If $\Gamma=\{1\}$ we simply have a fibration and the theorem
 reduces to the Atiyah-Singer family index theorem.
If $T=$ point then we have a Galois
covering and Connes' index theorem reduces to the Connes-Moscovici
higher index theorem.
 If $\dim T >0$ and $\Gamma\not= \{1\}$, then
Connes' index theorem  is a higher foliation index
theorem on the foliated manifold $(Y,\mathcal{F})$.

One particularly interesting higher index is the so-called Godbillon-Vey  index;
an alternative treatment of Connes' index formula in this particular case was given by Moriyoshi-Natsume
in \cite{MN}. Subsequently, Gorokhovsky and Lott \cite{Go-Lo}
gave a superconnection proof of Connes'
index theorem, including an explicit  formula for the Godbillon-Vey higher index. Leichtnam and Piazza
\cite{LPETALE}
extended Connes' index theorem to foliated bundles with boundary, using an extension of Melrose
$b$-calculus and the Gorokhovsky-Lott superconnection approach. Unfortunately, a key assumption
in \cite{LPETALE} is that the group $\Gamma$ be of polynomial growth. This excludes many interesting
examples and
higher indeces; in particular it excludes the possibility of proving a Atiyah-Patodi-Singer  formula
for the Godbillon-Vey higher index.

\bigskip
 {\it One primary objective of this article is to illustrate such a result. Complete proofs will appear in \cite{MP}.}

 \bigskip
In tackling the problem we  develop what we believe is a new approach to index theory on manifolds with boundary.
This can be summarized as follows. We  define a short
exact sequence of Banach algebras $$0\to \mathbf{ \mathfrak{J}} \to \mathbf{ \mathfrak{A}} \to
\mathbf{ \mathfrak{B}} \to 0$$ with $ \mathbf{ \mathfrak{J}}$ dense
and holomorphically closed in $C^* (X,\F)$ and with $ \mathbf{ \mathfrak{B}}$
depending only on boundary data. We prove that there are well defined
 Dirac index classes, denoted respectively
$\Ind (D,D^\pa)\in K_* (\mathbf{ \mathfrak{A}},\mathbf{ \mathfrak{B}}) $ and $\Ind(D)\in K_* (\mathbf{ \mathfrak{J}})$, and that these index classes correspond under excision; the relative class $\Ind (D,D^\pa)\in K_* (\mathbf{ \mathfrak{A}},\mathbf{ \mathfrak{B}}) $ is obtained by using the graph projection of $D$ and of $D^{\cyl}$ whereas the index class
$\Ind(D)\in K_* (\mathbf{ \mathfrak{J}})$ is obtained through the parametrix of $D$ and the associated remainders.
Next, for  (certain) cyclic $k$-cocycles
 defining a higher index in the closed case, let us name one of such cocycles
 $\tau_{k}$, we define
\begin{itemize}
\item a cyclic $k$-cocycle on $\mathbf{ \mathfrak{J}}$, still denoted $\tau_{k}$;
\item a eta cyclic cocycle $\sigma_{k+1}$ on $\mathbf{ \mathfrak{B}}$;   $\sigma_{k+1}$
(which thus depends solely on boundary data) is obtained by a sort of suspension  procedure involving $\tau_{k}$ and
a specific  1-cocycle $\sigma_1$ (Roe's 1-cocycle);
\item a {\it relative} cyclic $k$-cocyle $(\tau^r_{k}, \sigma_{k+1})$, with $\tau^r_k$ a cyclic cochain defined from
$\tau_k$ through a regularization \`a la Melrose.
\end{itemize}
The index formula in this context  is obtained by establishing
the equality
$$\langle \Ind (D),[\tau_k] \rangle\,=\,\langle \Ind (D,D^\pa), [\tau^r_{k}, \sigma_{k+1}] \rangle\,.$$
On the left hand side we have
 the absolute pairing, which is by definition the higher index.
On the right hand side we have the relative pairing; multiplying the operator by $s>0$,
using the definition of the relative pairing and taking the limit as $s\downarrow 0$ we obtain
the right hand side of the  Atiyah-Patodi-Singer index formula. The eta-correction term is obtained through the
eta cocycle $\sigma_{k+1}$.

\smallskip
We end this brief introduction by pointing out that relative pairings in K-theory and cyclic cohomology have already been successfully employed in the study of geometric and topological invariants of
elliptic operators. We particularly wish to mention here the paper by Lesch, Moscovici and Pflaum
\cite{lmpflaum}; in this interesting article the absolute and relative pairings  associated to a suitable
short exact sequence of algebras (this is a  short exact sequence of parameter
dependent pseudodifferential operators)
are used in order to define and study a generalization of the
divisor flow of Melrose on a closed compact manifold, see  \cite{melrose-eta} and also \cite{LPfla}.

\medskip
\noindent
{\it The results of  this article  first appeared in \cite{MP-arxiv}. }

\medskip
\noindent
{\bf Acknowledgements.} Most of this work has been done while the first author was
visiting Sapienza Universit\`a di Roma and the second author was visiting Keio University
and Nagoya University. We thank the Japan Society for the Promotion of Science (JSPS), 
Grants-in-Aid for Scientific Research, 
and the 21st century COE program at Keio 
for sponsoring most of these visits. Further financial support was provided by
{\it Istituto Nazionale di Alta Matematica}, through the GNSAGA, and the {\it Ministero dell'Istruzione,
dell'Universit\`a e della Ricerca (MIUR)} through the project {\it Spazi di moduli e Teoria di Lie}.
Part of this research was also carried out while the two authors were visiting jointly
 the Chern Institute in Tianjin and  the {\it Institut de Math\'ematiques de Jussieu} 
in Paris (\'Equipe
Alg\`ebres d'Op\'erateurs). 
We thank these institutions for  hospitality and financial support.
 Finally, it is a pleasure to thank Sergio Doplicher, Sacha Gorokhovsky, Eric Leichtnam, Henri Moscovici,
Toshikazu Natsume, John Roe and Xiang Tang for helpful 
discussions.

\section{Geometry of foliated bundles.}\label{sec:data}

\subsection{Manifolds with boundary}
Let now
$(M,g)$ be a riemannian manifold with boundary; the metric is assumed to be of product type in a collar
neighborhood $U\simeq [0,1]\times \partial M$ of
the boundary. Let $\tM$ be a Galois $\Gamma$-cover of $M$; we let $\tilde g$ be the lifted metric. We also consider  $\partial \tM$, the
boundary of $\tM$. Let $T$ be a smooth oriented compact manifold with an action
of $\Gamma$ by orientation-preserving diffeomorphisms. We assume that this action
is locally faithful, as in \cite{MN}, that is:
if $\gamma \in \Gamma$ acts as the identity map on an open set in $T$,
then $\gamma$ is the identity element in $\Gamma$.

Let $X_0=\tM\times_\Gamma T$; this is a manifold with boundary and the boundary $\partial X_0$ is equal to $\partial \tM \times_\Gamma T$.
$(X_0,\mathcal F_0)$ denotes  the associated foliated bundle. The leaves of
$(X_0,\mathcal F_0)$ are manifolds with boundary
endowed with a product-type metric. The boundary $\partial X_0$ inherits a  foliation $\mathcal F_\partial$.
The  cylinder $\RR\times \partial X_0$ also inherits a foliation $\mathcal F_{{\rm cyl}}$, obtained by crossing
the leaves of $\mathcal F_\partial$ with $\RR$. Similar considerations apply to the half cylinders
$(-\infty,0] \times \partial X_0$ and $\partial X_0 \times [0,+\infty)  \,. $

\subsection{Manifolds with cylindrical ends. Notation.}
We consider
$\tV:= \tM\cup_{\partial \tM} \left(   (-\infty,0] \times \partial\tM \right),$
endowed with the extended metric and the obviously extended $\Gamma$ action along the cylindrical
end. We consider $X:= \tV\times_\Gamma T$; this is a foliated bundle, with leaves manifolds with cylindrical ends.
We denote by $(X,\mathcal F)$ this foliation. Notice that $X=X_0 \cup_{\partial X_0} \left(   (-\infty,0] \times \partial X_0 \right)$;
moreover the foliation $\mathcal F$ is obtained by extending $\mathcal F_0$ on $X_0$ to $X$ via the product cylindrical
foliation $\mathcal F_{{\rm cyl}}$ on $(-\infty,0] \times \partial X_0$. We can write more suggestively:
$ (X,\mathcal F)= (X_0,\mathcal F_0)\cup_{(\partial X_0,\mathcal F_\partial)} \left(   ((-\infty,0] \times \partial X_0,  \mathcal F_{{\rm cyl}})\right)$.
For $\lambda>0$ we shall also consider the finite cyclinder $\tV_\lambda = \tM\cup_{\partial \tM} \left(   [-\lambda,0] \times \partial\tM \right)$
and the resulting foliated manifold $(X_\lambda,\mathcal F_\lambda)$.
Finally, with a small abuse, we  introduce the notation:
$\cyl (\pa X):= \RR\times \partial X_0$, $
\cyl^- (\pa X):=(-\infty,0] \times \partial X_0$ and $\cyl^+ (\pa X):=\partial X_0 \times [0,+\infty)$
The foliations induced on $\cyl (\pa X)$, $\cyl^\pm (\pa X)$ by $\F_{\pa}$ will be denoted by
$\F_{\cyl}$, $\F_{\cyl}^\pm$.

\subsection{Holonomy groupoid}
We consider the
groupoid $G:=(\tV\times\tV\times T)/\Gamma$
with $\Gamma$ acting diagonally; the source map and the range map
are defined  by
$s[y,y',\theta]=[y',\theta]$, $ r[y,y',\theta]=[y,\theta]$.
Since the action on $T$ is assumed to be locally faithful, we know that $(G,r,s)$
is isomorphic to the
holonomy groupoid of the foliation $ (X,\mathcal F)$.
In the sequel, we shall simply call
$(G,r,s)$ {\it the holonomy groupoid}.
If $E\to X$ is a hermitian vector bundle on $X$, with product structure along the cylindrical end, then
we can consider the bundle over $G$ equal to $(s^* E)^*\otimes r^*E$.

\section{Wiener-Hopf extensions}\label{sec:algebras}

\subsection{Foliation $C^*$-algebras}\label{subsec:cstar}
 We consider $C_c (X,\mathcal F):= C_c (G)$.
  $C_c (X,\mathcal F)$ can also be defined  as the space of $\Gamma$-invariant continuous functions on
$\tV\times\tV\times T$ with $\Gamma$-compact support. More generally we consider
$
C_c (X,\mathcal F;E) := C_c (G,(s^* E)^*\otimes r^*E )
$
with its well known  *-algebra structure given by convolution. We shall often omit the vector bundle $E$ from
the notation.

The foliation $C^*$-algebra
$C^*(X,\mathcal F;E)$ is defined by completion of $C_c (X,\mathcal F;E)$. See for example \cite{MN} where it is also proved that
$C^*(X,\mathcal F;E)$ is isomorphic to the $C^*$-algebra of compact operators of the Connes-Skandalis $C(T)\rtimes \Gamma$-Hilbert module
$\mathcal E$ (this is also described in \cite{MN}). Summarizing: $C^*(X,\mathcal F;E)\cong \KK (\mathcal E)\subset \mathcal L (\mathcal E)$.

\subsection{Foliation von Neumann algebras}

Consider the family of Hilbert spaces $\H:= (\H_\theta)_{\theta\in T}$, with
$\H_\theta := L^2 (\tilde{V}\times \{\theta\}, E_\theta)$.
Then $C_c (\tilde{V}\times T)$ is a continuous field inside $\H$, that is, a linear subspace in the
space of measurable sections of $\H$. Let $\End (\H)$ the space of measurable families
of bounded operators $T=(T_\theta)_{\theta\in T}$, where bounded means that each
$T_\theta$ is bounded on $\H_\theta$.
Then $\End (\H)$ is a $C^*$-algebra, in fact a von Neumann algebra, equipped with the norm
$$\| T \|:= {\rm ess.} \sup \{\|T_\theta\|\,;\theta\in T\}\,$$
with $\|T_\theta\|$ the operator norm. We also denote by $\End_\Gamma (\H)$ the $C^*$-subalgebra
of $\End (\H)$ consisting of $\Gamma$-equivariant measurable families  of operators. This is often
denoted $W^* (X,\F)$ and named the {\it foliation von Neumann algebra} associated to  $(X,\F)$.
We set $C^*_\Gamma (\H)$ the closure of $\Gamma$-equivariant families $T=(T_\theta)_{\theta\in T}
\in \End_\Gamma (\H)$ preserving the continuous field  $C_c (\tilde{V}\times T)$. In \cite{MN}, Section 2
it is proved that the foliation $C^*$-algebra $C^* (X,\F)$ is isomorphic to a $C^*$-subalgebra of
$C^*_\Gamma (\H)\subset \End_\Gamma (\H)$ \footnote{The
 $C^*$-algebra $C^*_\Gamma (\H)$ was denoted $\mathfrak{B}$ in \cite{MN}}. Notice, in particular,
that an element in $C^* (X,\F)$
can be considered as a $\Gamma$-equivariant
family of operators.

\subsection{Translation invariant operators}\label{subsec:translation}

Recall $\cyl (\pa X):=\RR\times \partial X_0\equiv (\RR\times\partial \tM)\times_\Gamma T$ with $\Gamma$ acting trivially
in the $\RR$-direction of $(\RR\times\partial \tM)$.
We consider the foliated cylinder $(\cyl (\pa X), \mathcal F_{{\rm cyl}})$ and its holonomy groupoid
$G_{{\rm cyl}}:=((\RR\times\partial \tM)\times (\RR\times\partial \tM)\times T)/\Gamma$
(source and range maps are clear). Let $\RR$ act trivially on $T$; then $(\RR\times\partial \tM)\times (\RR\times\partial \tM)\times T$
has a $\RR\times\Gamma$-action, with $\RR$ acting by translation on itself.
We consider the *-algebra
$B_c (\cyl (\pa X),\mathcal F_{{\rm cyl}})\equiv B_c$ defined as
\begin{equation}\label{eq:tras-algebra}
\{k\in C((\RR\times\partial \tM)\times (\RR\times\partial \tM) \times T); k
\text{ is } \RR\times \Gamma\text{-invariant}, k \text{ has } \RR\times\Gamma\text{-compact  support}\}
\end{equation}
The product is by convolution. An element $\ell$ in
$B_c$  defines
a $\Gamma$-equivariant
family $(\ell(\theta))_{\theta\in T}$ of translation-invariant operators. The completion of
 $ B_c$ with respect to the obvious $C^*$-norm (the sup over $\theta$
of the operator-$L^2$-norm of $\ell(\theta)$) gives us a $C^*$-algebra that will be denoted $B^{*} (\cyl (\pa X),\F_{{\cyl}})$
or more briefly $B^{*}$.

\subsection{Wiener-Hopf extensions}\label{subsec:extension}

Recall the Hilbert  $C(T)\rtimes\Gamma$-module $\mathcal{E}$ and the $C^*$-algebras $\KK(\mathcal{E})$ and $\mathcal{L}( \mathcal{E})$.
Since  the $C(T)\rtimes\Gamma$-compact operators $\KK(\mathcal{E})$ are an ideal in  $\mathcal{L}( \mathcal{E})$
 we have the classical  short exact sequence of $C^*$-algebras
$$0\to  \KK(\mathcal{E})\hookrightarrow \mathcal{L}( \mathcal{E}) \xrightarrow{\pi} \mathcal{Q} ( \mathcal{E})\rightarrow 0$$
with $\mathcal{Q} ( \mathcal{E})= \mathcal{L}( \mathcal{E})/\KK(\mathcal{E})$ the Calkin algebra.
Let $\chi_\RR^0:\RR\to \RR$ be the characteristic function of $(-\infty,0]$; let $\chi_\RR:\RR\to\RR$ be a smooth function with values in $[0,1]$ such that $\chi (t)=1$ for $t\leq -\epsilon$, $\chi (t)=0$ for $t\geq 0$.
Let $\chi^0$ and $\chi$ be the  functions induced by $\chi^0_\RR$ and $\chi_\RR$ on $X$.
Similarly, introduce   $\chi^0_{{\rm cyl}}$ and $\chi_{{\rm cyl}}$.
\begin{lemma}\label{lemma:split}
There exists a bounded linear map
\begin{equation}\label{eq:split-2}
s : B^* \to \mathcal{L}(\mathcal{E})
\end{equation}
extending $s_c : B_c \to
\mathcal{L} ( \mathcal{E})$,  $s_c (\ell) := \chi^0 \ell \chi^0 $.
Moreover, the composition $\rho =\pi s $ induces an {\bf injective}
$C^*$-homomorphism
\begin{equation}\label{eq:split}
\rho : B^* \to \mathcal{Q}(\mathcal{E}).
\end{equation}
\end{lemma}
We  consider  $\Im \rho$ as a $C^*$-subalgebra in $ \mathcal{Q} ( \mathcal{E})$ and
 identify it with  $B^* (\cyl (\pa X), \mathcal F_{{\rm cyl}})$ via $\rho$.
Set $$A^* (X;\mathcal F):= \pi^{-1} (\Im \rho) \;\;\text{ with }
\pi \;\text{ the projection }\;\;
\mathcal{L}( \mathcal{E}) \rightarrow \mathcal{Q} ( \mathcal{E}).$$
Recalling the identification $C^*(X,\mathcal F)=\KK(\mathcal{E})$,
we thus  obtain a short exact sequence of $C^*$-algebras:
 \begin{equation}\label{eq:short-cstar}
 0\rightarrow C^*(X,\mathcal F)\rightarrow
A^* (X;\mathcal F)\xrightarrow{\pi} B^* (\cyl (\pa X), \mathcal F_{{\rm cyl}}) \rightarrow 0
 \end{equation}
 where
  the quotient map is still denoted by $\pi$.
 Notice that \eqref{eq:short-cstar} splits as a short exact sequence of {\it Banach spaces},
since we can choose
$s: B^* (\cyl (\pa X), \mathcal F_{{\rm cyl}}) \to A^* (X;\mathcal F)$  the map in \eqref{eq:split-2} as a section.
So $$A^* (X;\mathcal F) \cong C^*(X,\mathcal F)\oplus s( B^* (\cyl (\pa X), \mathcal F_{{\rm cyl}}))$$ as Banach spaces.

  There is also a linear map $t: A^* (X,\mathcal F)\rightarrow C^*(X,\mathcal F)$ which is obtained as follows: if
 $k\in A^* (X;\mathcal F)$,  then  $k$
 is  uniquely expressed as
 $k= a + s(\ell)$ with $a\in C^*(X,\mathcal F)$ and
$\pi(k)=\ell\in B^* (\cyl (\pa X), \mathcal F_{{\rm cyl}}) $.
Thus,  $\pi(k)=[\chi^0 \ell \chi^0]\in\mathcal Q (\mathcal E )$ for  one (and only one)  $\ell\in B^* (\cyl (\pa X), \mathcal F_{{\rm cyl}})$
since $\rho$ is injective.
 We set
\begin{equation}\label{eq:t}
t(k):= k-s\pi (k) =k-\chi^0 \ell \chi^0
\end{equation}
Then  $t(k)\in C^*(X,\mathcal F)$.

\section{Relative pairings and the eta cocycle: the algebraic theory}\label{sect:algebraic}

\subsection{Introductory remarks}\label{subsection:intr-remarks}
On a {\it closed} foliated bundle $(Y,\mathcal{F})$, the Godbillon-Vey cyclic cocycle is initially defined
on the "small" algebra $\mathcal{A}_c \subset C^*(Y,\mathcal{F}) $  of  $\Gamma$-equivariant
smoothing operators of $\Gamma$-compact support (viz. $\mathcal{A}_c:=C^\infty_c (G, (s^*E)^*\otimes r^*E)$).
Since the index class defined using a pseudodifferential parametrix is already well defined in $K_* (\mathcal{A}_c)$, the pairing between the
the Godbillon-Vey cyclic cocycle and the index class is well-defined.

In a second stage, the cocycle is continuously extended to a dense holomorphically closed subalgebra $\mathfrak{A}\subset C^*(Y,\mathcal{F})$;
there are at least two reasons
for doing this. First, it is only by going to the $C^*$-algebraic index that the
well known properties for the signature  and the spin Dirac operator  of  a metric of
positive scalar curvature hold. The second reason
for this extension rests on the structure of the index class {\it which is employed in
the proof of  the higher index formula}, i.e. either the graph projection or
the Wassermann projection; in both cases $\mathcal{U}_c$ is too small to contain the  index class and one is therefore forced to find
an  intermediate subalgebra $\mathfrak{A}$,
$\mathcal{A}_c \subset\mathfrak{A} \subset C^*(Y,\mathcal{F})$;
{\it $\mathfrak{A}$  is big enough for  the two particular index classes to belong to it but small enough for the Godbillon-Vey cyclic cocycle to extend;
moreover, being dense and holomorphically closed it has the  same $K$-theory as $C^*(Y,\mathcal{F})$.}

Let now $(X,\mathcal{F})$ be a foliated bundle with cylindrical ends; in  this section we shall select "small" subalgebras
$$J_c\subset  C^*(X,\mathcal F)\,,\quad A_c\subset  A^*(X,\mathcal F)\,,\quad B_c\subset B^* (\cyl (\pa X), \mathcal F_{{\rm cyl}})\,,$$
with $J_c$ an ideal in $A_c$,
so that there is a short exact sequence
$0\rightarrow J_c\hookrightarrow A_c \xrightarrow{\pi_c} B_c \rightarrow 0$
which is a subsequence of
 $0\rightarrow C^*(X,\mathcal F)\hookrightarrow A^* (X;\mathcal F)\xrightarrow{\pi} B^* (\cyl (\pa X), \mathcal F_{{\rm cyl}})$ $\rightarrow 0$.
We shall then proceed to define the relevant cyclic cocycles, relative and absolute, and study, algebraically,
their main properties. As in the closed case, we shall eventually need to find an intermediate short exact sequence,
sitting between the two, 
$0\rightarrow \mathfrak{J} \hookrightarrow \mathfrak{A} \rightarrow  \mathfrak{B} \rightarrow 0$,
with constituents big enough for the the two index classes we shall define to belong to them but small enough for
the cyclic cocycles (relative and absolute) to extend; this is quite a delicate point and it will be explained  in  Section
\ref{sec:shatten}.
We anticipate that in contrast with the closed case the ideal $J_c$ in the small subsequence will be too small even for
the index class defined by a pseudodifferential parametrix. This has to do with the non-locality of the parametrix on manifolds with boundary;
it is a phenomenon that was explained in detail in \cite{LPETALE}.

\subsection{Small dense subalgebras}\label{subsect:small-dense}

Define $J_c:= C^\infty_c(X,\mathcal{F})$; see subsection \ref{subsec:cstar}. Redefine $B_c$ as
$$ \{k\in C^\infty((\RR\times\partial \tM)\times (\RR\times\partial \tM) \times T); k
\text{ is } \RR\times \Gamma\text{-invariant}, k \text{ has } \RR\times\Gamma\text{-compact  support}\}
$$
see subsection \ref{subsec:translation}.
We now define $A_c$; consider the functions $\chi^\lambda$, $\chi^\lambda_{{\rm cyl}}$ induced on $X$ and $\cyl (\pa X)$
by the real function $\chi^\RR_{(-\infty,-\lambda] }$. We shall say that
 $k$ is in $A_c$ if it is a smooth function on  $\tV\times\tV\times T$ which is
$\Gamma$-invariant and
for which there  exists $\lambda\equiv \lambda(k)>0$, such that
\begin{itemize}
\item $k-\chi^\lambda k \chi^\lambda$ is of $\Gamma$-compact support
\item there exists $\ell \in B_c$ such that $\chi^\lambda k \chi^\lambda =  \chi^\lambda_{{\rm cyl}} \ell \chi^\lambda_{{\rm cyl}}
$ on $((-\infty,-\lambda]\times \partial \tM) \times ((-\infty,-\lambda]\times \partial \tM) \times T$
\end{itemize}

\begin{lemma}\label{lemma:ac}
$A_c$ is a *-subalgebra of $A^*(X,\mathcal{F})$.
Let $\pi_c:=\pi |_{A_c}$; there is a short exact sequence of *-algebras
\begin{equation}\label{short-compact}
0\rightarrow J_c\hookrightarrow A_c\xrightarrow{\pi_c} B_c \rightarrow 0\,.
\end{equation}
\end{lemma}

\begin{remark}
Notice that the image of $A_c$ through $t|_{A_c}$ is not contained in $J_c$ since
$\chi^0$ is not even continuous. Similarly, the image of $B_c$ through $ s|_{B_c}$
is not contained in $A_c$.
\end{remark}

\subsection{Relative cyclic cocycles}

Let $A$ be a
$k-$algebra over $k=\CC$. Recall the
cyclic cohomology groups $HC^*(A)$  \cite{Co}.
Given a second algebra $B$ together with a surjective homomorphism
$\pi: A\to B$, one can define the relative cyclic complex
$C^n_\lambda(A,B):=\{(\tau,\sigma)\,:\, \tau\in C^n_\lambda (A), \sigma\in C^{n+1}_\lambda (B)\}$
with  coboundary map given by
$(\tau,\sigma)\longrightarrow (\pi^*\sigma - b\tau, b\sigma)\,.$
A relative cochain $(\tau,\sigma)$
is thus a cocycle if $b\tau=\pi^* \sigma$ and $b\sigma=0$.
One obtains in this way the relative cyclic cohomology groups $HC^* (A,B)$.
If $A$ and $B$ are Fr\'echet algebra, then we can also define the topological (relative) cyclic
cohomology groups. More detailed information are given,
for example,
in \cite{lmpflaum}.

\subsection{Roe's 1-cocycle}

In this subsection, and in the next two, we study a particular but important example.
We assume that $T$ is a point and that $\Gamma=\{1\}$, so that we are really considering
a compact manifold $X_0$ with boundary $\pa X_0$ and associated manifold with  cylindrical ends
$X$; we keep denoting the cylinder $\RR\times \pa X_0$ by $\cyl (\pa X)$ (thus, as before,
we don't write the subscript $0$).
 The algebras appearing in the short exact sequence \eqref{short-compact}
 are now given by
 $J_c=C^\infty_c (X\times X)$, $$B_c=\{k\in C^\infty((\RR\times \pa X_0)\times (\RR\times \pa X_0)); k \text{ is $\RR$-invariant}, k
 \text{ has compact $\RR$-support}\}\,.$$
Finally, a smooth function $k$ on $X\times X$ is in  $A_c$ if there exists
a $\lambda\equiv \lambda(k)>0$ such that \\(i) $k-\chi^\lambda k \chi^\lambda$ is
of compact support on $X\times X$; \\(ii) $\exists$ $\ell\in B_c$ such that  $\chi^\lambda k \chi^\lambda=
\chi^\lambda_{{\rm cyl}} \ell  \chi^\lambda_{{\rm cyl}}$ on $((-\infty,-\lambda]\times \pa X_0)\times (-\infty,-\lambda]\times \pa X_0)$\,.\\
For such a $k\in A_c$ we define $\pi_c (k)=\ell$ and  we have the short exact sequence of $*$-algebras
$0\rightarrow J_c\hookrightarrow A_c\xrightarrow{\pi_c} B_c \rightarrow 0\,.$
Incidentally, in  the Wiener-Hopf short exact sequence  \eqref{eq:short-cstar}, which now reads as
$0\rightarrow C^*(X)\rightarrow
A^* (X)\xrightarrow{\pi} B^* (\cyl (\pa X)) \rightarrow 0$, the left term $C^*(X)$ is clearly given by the
compact operators on $L^2 (X)$.

We shall define below a 0-{\it relative} cyclic cocycle associated to the homomorphism $\pi_c: A_c \to B_c$.
To this end we start by
defining  a cyclic 1-cocycle $\sigma_1$ for the algebra $B_c$; this is directly inspired from work
of  John Roe
(indeed, a similarly defined 1-cocycle  plays a fundamental  role in his index
theorem on partioned manifolds \cite{roe-partitioned}).

Consider the characteristic function $\chi^\lambda_{{\rm cyl}}$, $\lambda>0$,
 induced on the cylinder by the real function $\chi^\RR_{(-\infty,-\lambda]}$.
For notational convenience, unless absolutely necessary,
 we shall not distinguish between $\chi^\lambda_{{\rm cyl}}$ on the cylinder $\cyl (\pa X)$ and $\chi^\lambda$ on the manifold
with cylindrical ends $X$.

We define $\sigma_1^\lambda: B^{\RR}_c\times B_c \to \mathbb{C}$ as
\begin{equation}\label{roe-1}
\sigma^\lambda_1 (\ell_0,\ell_1):= \Tr (\ell_0 [\chi^\lambda\,,\,\ell_1])\,.
\end{equation}
One can check that the operators $[\chi^\lambda\,,\,\ell_0]$ and
$\ell_0 [\chi^\lambda\,,\,\ell_1]$ are trace class $\forall \ell_0, \ell_1 \in B_c$
(and
 $\Tr [\chi^\lambda\,,\,\ell_0]=0$). In particular
$\sigma^\lambda_1 (\ell_0,\ell_1)$ is well defined.

\begin{proposition}\label{prop:lambda-indipendent}
The value $\Tr (\ell_0 [\chi^\lambda\,,\,\ell_1])$ is independent of $\lambda$
and will simply be denoted by $\sigma_1 (\ell_0,\ell_1)$.
The functional $\sigma_1 :B_c\times B_c\to \CC$ is a cyclic 1-cocycle.
\end{proposition}

\subsection{Melrose' regularized integral}

Recall that our immediate goal is to define a {\it relative} cyclic 0-cocycle for the homomorphism
$\pi_c : A_c \to B_c$ appearing in the short exact sequence of the previous section. Having defined
a 1-cocycle $\sigma_1$ on $B_c$ we now need to define a 0-cochain on $A_c$. Our definition will be
a simple adaptation of the definition of the $b$-trace in Melrose' $b$-calculus \cite{Melrose}
(but since our algebra $A_c$ is very small, we can give a somewhat simplified treatment). Recall that
for $\lambda>0$ we are denoting by $X_\lambda$ the compact manifold obtained attaching $[-\lambda,0]\times \pa X_0$
to our manifold with boundary $X_0$.

So, let $k\in A_c$ with $\pi_c (k)=\ell\in B_c$. Since $\ell$ is $\RR$-invariant on the cylinder $\cyl (\pa X)=\RR\times \pa X_0$
we can write $\ell(y,y^\prime,s)$ with  $y,y^\prime\in \pa X_0$, $s\in \RR$.
Set
\begin{equation}\label{tau0}
\tau_0^r (k) := \lim_{\lambda\to +\infty} \left( \int_{X_\lambda} k(x,x) {\rm dvol}_g -\lambda \int_{\pa X_0} \ell(y,y,0) {\rm dvol}_{g_\partial}
\right)
\end{equation}
where the superscript $r$ stands for {\it regularized}. (The $b$-superscript would be of course more appropiate;
unfortunately it gets confused with the $b$ operator in cyclic cohomology.)
It is elementary to see that the limit exists; in fact, because of the very particular definition of
$A_c$ the function
$\varphi (\lambda):= \int_{X_\lambda} k(x,x) {\rm dvol}_g -\lambda \int_{\pa X_0} \ell(y,y,0) {\rm dvol}_{g_\partial}$
becomes constant for large values of $\lambda$. The proof is elementary.
$\tau^r_0$ defines a 0-cochain on $A_c$.

\begin{remark}
Notice that \eqref{tau0}  is nothing but Melrose' regularized integral \cite{Melrose}, in the cylindrical language,
for the restriction of $k$ to the diagonal of $X\times X$.
\end{remark}

We shall also need  the following

\begin{lemma}\label{lemma:regularized}
If $k\in A_c$ then $t(k)$, which is a priori a compact operator, is in fact trace class
and
$\tau^r_0 (k)= \Tr (t (k))\,.$
\end{lemma}
We remark once again that $t(k)$ is not an element in $J_c$.

\subsection{The regularized integral  and Roe's 1-cocycle define a relative 0-cocycle}

We finally consider the relative 0-cochain $(\tau_0^r,\sigma_1)$ for the pair $A_c\xrightarrow{\pi_c} B_c$.
\begin{proposition}\label{prop:c-relative-cocycle}
The relative 0-cochain $(\tau_0^r,\sigma_1)$ is a relative 0-cocycle. It thus defines an element $[(\tau_0^r,\sigma_1)]$
 in the relative group $HC^0 (A_c,B_c)$.
\end{proposition}
There are several proofs of this Proposition; we have stated that $\sigma_1$
is a cocycle and what needs to be proved now is that $b \tau^r_0= (\pi_c)^* \sigma_1$.
One proof of this equality employs Lemma \ref{lemma:regularized}; another one
use the Hilbert transform and Melrose' formula for the $b$-trace of a commutator \cite{Melrose}, see the next Subsection.

\subsection{Melrose' 1-cocycle and the relative cocycle condition via the $b$-trace formula}
As we have anticipated in the previous subsection, the equation $b\tau_0^r= \pi^*_c \sigma_1$ is nothing
but a compact way of rewriting Melrose' formula for the $b$-trace of a commutator. We wish to
explain this point here.
Following now the notations of the $b$-calculus, we consider the sligthy larger  algebras
$$
A^b_c:= \Psi^{-\infty}_{b,c}(X,E)\,,\quad B^b_c:= \Psi^{-\infty}_{b,I,c}(\overline{N_+\partial X},E|_{\pa})\,,\quad
J^b_c := \rho_{{\rm ff}}
\Psi^{-\infty}_{b,c}(X,E)$$
and  $
0\,\longrightarrow \,
J^b_c\,
\longrightarrow
A^b_c\,
\xrightarrow{ \pi^b_c}
B^b_c\,
\longrightarrow \,0$, with $ \pi^b_c $ equal to Melrose' indicial operator $I(\cdot)$. Let $\tau_0^r$ be equal to the $b$-Trace:
$\tau_0^r:={}^b \Tr$.
 Observe that $\sigma_1$ also defines a 1-cocyle on $B_c^b$.
 We can thus consider the relative 0-cochain $(\tau_0^r,\sigma_1)$ for the homomorphism
$A^b_c\,
\xrightarrow{I(\cdot)}
B^b_c$; in order to prove that this is  a relative 0-cocycle it remains to
 to show that $b\tau_0^r (k,k^\prime)= \sigma_1 (I(k),I(k^\prime))$, i.e.
 \begin{equation}\label{1cocycle-is-b}
 {}^b \Tr [k,k^\prime]= \Tr (I(k)[\chi^0,I(k^\prime)])
 \end{equation}

Recall here that Melrose' formula for the $b$-trace of a commutator is
 \begin{equation}\label{melrose}
 {}^b \Tr [k,k^\prime]= \frac{i}{2\pi}\int_\RR \Tr_{\pa X} \left( \partial_\mu I(k,\mu) \circ I (k^\prime,\mu) \right) d\mu
 \end{equation}
 with
 $\CC\ni z\to I(k,z)$ denoting the indicial family of the operator $k$, i.e. the Fourier transform
of the indicial operator $I(k)$.

Inspired by the right hand side of \eqref{melrose} we consider an arbitrary  compact manifold
 $Y$, the algebra $B^b_c (\cyl (Y))$ and the 1-cocycle
\begin{equation}\label{melrose-1-cocycle}
\mathfrak{s}_1 (\ell,\ell^\prime):=  \frac{i}{2\pi}\int_\RR \Tr_Y \left( \partial_\mu \hat{\ell}(\mu) \circ \hat{\ell}^\prime (\mu) \right) d\mu
\end{equation}
That this is a cyclic 1-cocyle follows by elementary arguments.
Formula \eqref{melrose-1-cocycle}
defines what should be called Melrose' 1-cocycle

\begin{proposition}\label{prop:hilbert-transf}
Roe's 1-cocycle $\sigma_1$ and Melrose 1-cocycle $\mathfrak{s}_1 $
coincide:
\begin{equation}\label{equality-with-rbm-bis}
\sigma_1 (\ell,\ell^\prime):= \Tr (\ell[\chi^0,\ell^\prime]) = \frac{i}{2\pi}\int_\RR \Tr_Y \left( \partial_\mu \hat{\ell}(\mu) \circ \hat{\ell}^\prime (\mu) \right) d\mu
 =:  \mathfrak{s}_1 (\ell,\ell^\prime)
\end{equation}
\end{proposition}
Proposition \ref{prop:hilbert-transf} and Melrose' formula imply at once the relative 0-cocyle
condition for $(\tau_0^r,\sigma_1)$:
indeed using first Proposition \ref{prop:hilbert-transf} and then Melrose' formula we get:
\begin{align*}
\sigma_1 (I(k),I(k^\prime))&:=\Tr (I(k)[\chi^0,I(k^\prime)])= \frac{i}{2\pi}\int_\RR \Tr_{\pa X} \left( \partial_\mu I(k,\mu) \circ I (k^\prime,\mu) \right) d\mu
\\
&= {}^b \Tr [k,k^\prime]= b\tau_0^r (k,k^\prime)\,.
\end{align*}
Thus $I^* (\sigma_1)=b\tau_0^r$ as required.

\noindent
{\bf Conclusions.} We have seen the following:
\begin{itemize}
\item
the right hand side of Melrose' formula defines a 1-cocyle $\mathfrak{s}_1$ on $B_c (\cyl (Y))$,
$Y$ any closed compact manifold;
\item  Melrose 1-cocyle $\mathfrak{s}_1$ equals Roe's 1-cocyle $\sigma_1$
\item
Melrose' formula itself can be interpreted
as a {\it relative} 0-cocyle condition for the 0-cochain $(\tau^r_0,\mathfrak{s}_1)\equiv (\tau^r_0,\sigma_1)$.
\end{itemize}

\subsection{Philosophical remarks}\label{subsect:strategy}

We wish to recollect the information obtained in the last three subsections and start to explain
our approach to Atiyah-Patodi-Singer higher index theory.

On a closed compact  orientable riemannian smooth manifold $Y$ let us consider the algebra of smoothing operators
$J_c (Y):= C^\infty (Y\times Y)$.
Then the functional $J_c (Y) \ni k\rightarrow \int_Y k|_\Delta {\rm dvol}$
defines a 0-cocycle $\tau_0$ on $J_c (Y)$; indeed by Lidski's theorem the functional is nothing but
the functional analytic trace of the integral operator corresponding to $k$ and we know that
the trace vanishes on commutators; in formulae, $b\tau_0 =0$. The 0-cocycle $\tau_0$
plays a fundamental role in the proof of the Atiyah-Singer index theorem, but we leave this
aside for the time being.

Let now $X$ be a smooth orientable manifold with cyclindrical ends, obtained from a manifold
with boundary $X_0$; let $\cyl (\pa X)=\RR\times \pa X_0$.
 We have then defined algebras $A_c (X)$, $B_c (\cyl (\pa X))$ and $J_c (X)$ fitting into a short exact sequence
 $0\to J_c (X)\to A_c (X)\xrightarrow{\pi_c}B_c (\cyl (\pa X))\to 0$.

 Corresponding to the 0-cocycle $\tau_0$ in the closed
 case we can define two important 0-cocycles  on a manifold with cyclindrical ends
$X$:

\begin{itemize}
\item We can consider $\tau_0$ on $J_c (X)=C^\infty_c (X\times X)$; this is well defined
and does define a 0-cocycle . We shall refer to $\tau_0$ on
$J_c (X)$ as an {\it absolute} 0-cocycle.
\item Starting with the absolute 0-cocycle $\tau_0$ on $J_c (X)$ we
 define a {\it relative} 0-cocycle $(\tau_0^r,\sigma_1)$ for $A_c (X)\xrightarrow{\pi_c}B_c (\cyl (\pa X))$.
  The relative 0-cocycle   $(\tau_0^r,\sigma_1)$ is obtained
 through the following two fundamental steps.
 \begin{enumerate}
 \item We  define a 0-{\it cochain} $\tau_0^r$ on $A_c (X)$ by replacing the integral with Melrose' regularized
 integral.
 \item We define a 1-cocycle $\sigma_1$ on $B_c (\cyl (\pa X))$ by
 taking
a {\it suspension}
 of $\tau_0$ through the linear map $\delta(\ell):= [\chi^0,\ell]$. In other words,   $\sigma_1 (\ell_0,\ell_1)$
  is obtained from $\tau_0\equiv \Tr$ by
 considering $(\ell_0,\ell_1)\rightarrow \tau_0 (\ell_0 [\chi^0,\ell_1])\equiv \tau_0 (\ell_0 \delta(\ell_1))$.
\end{enumerate}
  \end{itemize}

\begin{definition}
We shall also call   Roe's 1-cocycle $\sigma_1$ {\it the  eta 1-cocycle corresponding to the absolute 0-cocycle $\tau_0$}.
\end{definition}

In order to justify the wording of this definition we need to show that all this has something to
do with the eta invariant and its role in the Atiyah-Patodi-Singer index formula. This will be  explained
in Section \ref{sec:index} and Section \ref{sect:index-theorems}.

\subsection{The absolute Godbillon-Vey 2-cyclic cocycle $\tau_{GV}$}\label{subsect:absolute-taugv}

Let $(Y,\F)$, $Y=\tN\times_\Gamma T$, be a foliated bundle {\it without} boundary.
We assume that $T=S^1$.
Let $E\to Y$ an hermitian complex vector bundle
on $Y$.  Let $\widehat{E}$ be the $\Gamma$-equivariant lift of $E$ to
$\tN\times T$.
Let $(G,s:G\to Y,r:G\to Y)$ be the holonomy groupoid
associated to $Y$, $G=(\tN\times\tN\times T)/\Gamma$. Consider again  the convolution algebra
$C^\infty_c (G, (s^*E)^*\otimes r^*E)$,
of equivariant smoothing
families with $\Gamma$-compact support.
The notation
$\Psi^{-\infty}_c (G,E)$ is also employed.
On  $\Psi^{-\infty}_c (G,E)\equiv C^\infty_c (G, (s^*E)^*\otimes r^*E)$ we can define a remarkable 2-cocycle, denoted $\tau_{GV}$ and
known as the Godbillon-Vey cyclic cocycle.
First of all, recall that there is a
weight $\omega_{\Gamma}$
defined on the
algebra $\Psi^{-\infty}_c (G;E)$,
\begin{equation}\label{weight-again}
\omega_{\Gamma}(k)=
 \int_{Y(\Gamma)} \mathrm{Tr}_{(\tilde{n},\theta)} k(\tilde{n},\tilde{n},\theta)d\tilde{n}\,d\theta\,.
\end{equation}
In this formula $Y(\Gamma)$ is the fundamental domain in $\tN\times T$ for the free diagonal action
of $\Gamma$  on $\tN\times T$ and  we have restricted the kernel $k$ to $\Delta_{\tN}\times T\subset
\tN\times\tN\times T$,  $\Delta_{\tN}$ denoting the diagonal
set in $\tN\times \tN$, $\Delta_{\tN}\times T\equiv \tN\times T$;
$\mathrm{Tr}_{(\tilde{n},\theta)}$ denotes the trace on $\End (\widehat{E}_{(\tilde{n},\theta)})$.
(If the measure on $T$ is $\Gamma$-invariant, then this weight is
a trace; however, we don't want to make this assumption here.)
Recall then the bundle $\widehat{E}^\prime$ on $Y\times T$: this is the same vector bundle as
$\widehat{E}$ but with a different $\Gamma$-action. See \cite{MN} for details. There is a natural
identification $\Psi^{-\infty}_c (G;E)\equiv \Psi^{-\infty}_c (G;E^\prime)$.
We shall consider the linear space $\Psi^{-\infty}_c (G;E,E^\prime)$; using the
above identification we can give  $\Psi^{-\infty}_c (G;E,E^\prime)$  a natural bimodule
structure over $\Psi^{-\infty}_c (G;E)$.
We shall be interested in the linear functional \footnote{This will not be a weight, given
that on a bimodule there is no notion of positive element} defined on the bimodule
 $\Psi^{-\infty}_c (G;E,E^\prime)$
by the analogue of \eqref{weight-again}. To be quite explicit
\begin{equation}\label{weight-bimodule}
\omega_{\Gamma}(k)=
 \int_{Y(\Gamma)} \mathrm{Tr}_{(\tilde{n},\theta)} k(\tilde{n},\tilde{n},\theta)d\tilde{n}\,d\theta\,,\quad k\in \Psi^{-\infty}_c (G;E,E^\prime)
\end{equation}
where we now identify $\widehat{E}_{(\tilde{n},\theta)}$ and $\widehat{E}_{(\tilde{n},\theta)}^\prime$
given that thet are identical vector spaces (it is only the $\Gamma$-actions that are different).
We call \eqref{weight-bimodule} the  {\it bimodule trace}. (This name come from the following fundamental
property: if
$k\in \Psi^{-\infty}_c (G;E,E^\prime)$, $k^\prime \in \Psi^{-\infty}_c (G;E)\equiv \Psi^{-\infty}_c (G;E^\prime)$ then $\omega_{\Gamma}(k k^\prime)=\omega_{\Gamma}( k^\prime k)$.)

Recall now  the  two derivations $\delta_2:= [\phi,\;\;]$ and $\delta_1:=[\dot{\phi},\;\;]$ coming from the modular
automorphism group described in \cite{MN}.
More precisely,  we have a derivation $\delta_2$ and  a bimodule derivation $\delta_1$,
\begin{equation}\label{deltas}
\delta_2:  \Psi^{-\infty}_c (G,E)\to \Psi^{-\infty}_c (G;E)\,,\quad \delta_1 : \Psi^{-\infty}_c (G,E)\to \Psi^{-\infty}_c (G;E,E^\prime)\,,
\end{equation}

\begin{definition}\label{def:gv-cocycle}
With $1=\dim T$,
the Godbillon-Vey cyclic $2$-cocycle
on  \\$C^\infty_c (G, (s^*E)^*\otimes r^*E)$  is defined to be
\begin{equation}\label{taugv}
\tau_{GV} (a_0, a_1, a_2)= \frac{1}{2 !}\sum_{\alpha\in \mathfrak{S}_2} {\rm sign}(\alpha)\,
\omega_\Gamma (a_0 \;\delta_{\alpha (1)} a_1 \;\delta_{\alpha (2)} a_2)
\end{equation}
with $\omega_\Gamma$ the bimodule trace in   \eqref{weight-bimodule}.

\end{definition}

The fact that this 3-linear functional is indeed a cyclic 2-cocycle is proved in \cite{MN}.
We now go back to a foliated bundle $(X,\F)$ with cylindrical ends, with $X:= \tM\times_\Gamma T$, as in Section
\ref{sec:data}. We consider the small subalgebras introduced in Subsection \ref{subsect:small-dense}.
The weight $\omega_\Gamma$ is still well defined on $J_c (X,\F)$;
the 2-cocycle $\tau_{GV}$ can thus be defined on $J_c (X,\F)$, giving us the {\it absolute} Godbillon-Vey
cyclic cocycle.

\subsection{The eta 3-coycle $\sigma_{GV}$ corresponding to $\tau_{GV}$}\label{subsection:sigma3}

Now we apply the general philosophy explained at the end of the previous Section.
Let $\chi^0$ be the usual characteristic function of $(-\infty,0]\times \pa X_0$ in $\cyl (\pa X)=\RR\times \pa X_0$.
Write $\cyl (\pa X)= (\RR\times\pa \tM)\times_\Gamma T$ with $\Gamma$ acting trivially on the $\RR$ factor.
Let $\cyl (\Gamma)$ be a fundamental domain for the action of $\Gamma$ on $(\RR\times\pa \tM)\times T$;
finally,
let $\omega_\Gamma^{\,{\rm cyl}}$ be the corresponding weight. We keep
denoting this weight by $\omega_\Gamma$.
Recall the derivation $\delta (\ell): =[\chi^0,\ell]$;
recall that we passed from the absolute 0-cocycle $\tau_0\equiv \Tr$
to the 1-eta cocycle on the cylindrical algebra $B_c$ by considering $(\ell_0, \ell_1)\rightarrow \tau_0 (\ell_0 \delta (\ell_1))$.
We referred to this operation as a {\it suspension}.

We are thus led  to {\it suspend}
definition \ref{def:gv-cocycle}, thus defining the following 4-linear functional on the algebra $B_c$.

\begin{definition}\label{def:sigma3}
 The eta cochain $\sigma_{GV}$ associated to the absolute
Godbillon-Vey 2-cocycle $\tau_{GV} (a_0, a_1, a_2)$
is by definition
\begin{equation}\label{sigma3}
\sigma_{GV} (\ell_0, \ell_1, \ell_2,\ell_3)= \frac{1}{3 !}\,\sum_{\alpha\in \mathfrak{S}_3} \,{\rm sign}(\alpha)\,
\omega_\Gamma (\ell_0 \;\delta_{\alpha (1)} \ell_1\; \delta_{\alpha (2)} \ell_2 \;\delta_{\alpha(3)} \ell_3)
\end{equation}
with $\delta_3 (\ell):= [\chi^0,\ell]$.
The eta cochain is a 4-linear functional on $B_c (\cyl (\pa X),\F_{{\rm cyl}}))$
\end{definition}

In fact, we can define, as we did for $\sigma_1$, the 3-cochain $\sigma_{GV}^\lambda$  by
employing the characteristic function $\chi^\lambda$. However, one checks easily that the value of
$\sigma_{GV}^\lambda$ does not depend on $\lambda$.
One can prove that this definition is well posed, namely that
each term
$(\ell_0 \;\delta_{\alpha (1)} (\ell_1)\; \delta_{\alpha (2)} (\ell_2) \;\delta_{\alpha(3)} (\ell_3))$
is of finite weight.
We then have the important

\begin{proposition}\label{prop:sigma3-cyclic}
The eta functional $\sigma_{GV}$ is cyclic and it is a 3-cocycle: $b \,\sigma_{GV}=0$.
\end{proposition}

\subsection{The relative Godbillon-Vey cyclic cocycle $(\tau^r_{GV},\sigma_{GV})$}

We now apply our strategy as in Subsection \ref{subsect:strategy}. Thus starting with the absolute cyclic cocycle
$\tau_{GV}$ on $J_c (X,\F)$
 we first consider the 3-linear functional  on $A_c (X,\F)$ given by
$
 \psi_{GV}^r (k_0,k_1,k_2) :=  \frac{1}{2 !}\, \sum_{\alpha\in \mathfrak{S}_2} {\rm sign}(\alpha)\,
\omega_\Gamma^r (k_0 \;\delta_{\alpha (1)} k_1 \;\delta_{\alpha (2)} k_2)
$
with $\omega^r_{\Gamma}$ the regularized weight corresponding to $\omega_{\Gamma}$\\
 Next we consider the {\it cyclic} cochain associated to  $\psi_{GV}^r $:
\begin{equation}\label{regularized-gv}
\tau_{GV}^r (k_0,k_1,k_2):=   \frac{1}{3} \left( \psi_{GV}^r (k_0,k_1,k_2) +  \psi_{GV}^r (k_1,k_2,k_0)
+  \psi_{GV}^r (k_2,k_0,k_1) \right) \,.
\end{equation}

The next Proposition  is crucial:

\begin{proposition}\label{prop:gv-cocycle}
The relative cyclic cochain $(\tau_{GV}^r,\sigma_{GV})\in C^2_\lambda (A_c,B_c)$ is a relative
2-cocycle: thus $b\sigma_{GV}=0$ (which we already know) and $b\tau_{GV}^r= (\pi_c)^* \sigma_{GV}$.
\end{proposition}

  For later use we also state the analogue of
 Lemma \ref{lemma:regularized}:
 \begin{proposition}\label{prop:regularized-via-t}
 Let $t: A^* (X,\F)\to C^* (X,\F)$ be  the section introduced in Subsection \ref{subsec:extension}.
 If $k\in A_c\subset A^* (X,\F)$ then $t(k)$ has finite weight.
Moreover, for  the regularized weight $\omega_\Gamma^r: A_c \to \CC$
 we have
 \begin{equation}\label{regularized-via-t}
\omega^r_\Gamma = \omega_\Gamma \circ t
\end{equation}

 \end{proposition}

 \section{Smooth subalgebras}\label{sec:shatten}
In this section we select  important
 subsequences of $ 0\rightarrow C^*(X,\mathcal F)\rightarrow A^* (X;\mathcal F)\rightarrow B^* (\cyl (\pa X), \mathcal F_{{\rm cyl}}) \rightarrow 0$.

\subsection{Shatten ideals}\label{subsect:shatten-ideals}
Let $\chi_\Gamma$ be a characteristic function for a fundamental domain of $\Gamma\to \tM\to M$.
Consider $\Psi^{-\infty}_c (G;E)=:J_c (X\,\F)\equiv J_c$. We shall often omit the bundle $E$ from the
notation.

\begin{definition}\label{def:shatten}
Let $k\in J_c$ be positive and self-adjoint. The Schatten norm $||k||_m$ of $k$ is defined as
\begin{equation}\label{shatten-norm}
(||k||_m)^m := 
\sup_{\theta\in T} || \chi_\Gamma \, (k(\theta))^m \chi_\Gamma ||_1
\end{equation}
with the $||\,\,||_1$ denoting the usual trace-norm on the Hilbert space $\mathcal{H}_\theta=L^2(\tV\times\{\theta\}$.
Equivalently
\begin{equation}\label{shatten-norm-HS}
(||k||_m)^m = \sup_{\theta\in T} || \chi_\Gamma \, (k(\theta))^{m/2}  ||^2_{HS}\,.
\end{equation}
 with
 $|| \,\,\,||_{HS}$ denoting the usual Hilbert-Schmidt norm.
 In general, we set $||k||_m:= ||\,\,(k k^*)^{1/2}\,\,||_m$.
The Schatten norm of $k\in J_c$ is easily seen to be finite for any $m\geq 1$; we define $\mathcal{I}_m (X,\F)\equiv
\mathcal{I}_m$ as the completion of $J_c$
with respect to $||\;\;\;||_m$
\end{definition}
One can prove that $\mathcal{I}_m$ is a Banach algebra and  an ideal inside $C^*(X,\mathcal F)$.
Moreover:
\begin{proposition}\label{prop:extension-of-weight}
The weight $\omega_\Gamma$  extends continuously from $J_c\equiv C^\infty_c (G)$ to
 $\mathcal{I}_1$.
 \end{proposition}

We shall now introduce the subalgebra of $C^* (X,\F)$  that will be used in the proof of our index theorem.
Consider on the cylinder $\RR\times Y$ (with cylindrical variable $s$) the functions
\begin{equation}\label{functions-for-ideal}
f_{\cyl} (s,y):= \sqrt{1+s^2} \quad \quad g_{\cyl}(s,y)=1+ s^2\,.
\end{equation}
We denote by $f$ and $g$  smooth functions on $X$ equal to  $f_{\cyl}$ and $g_{\cyl}$ on
the open subset $(-\infty,0)\times Y$;
$f$ and $g$ are well defined up to a compactly supported function. We set
\begin{equation}\label{new-j}
\J_m (X,\F):= \{k\in\I_m \;|\: gk \;\text{and}\;kg \;\text{are bounded}\}
\end{equation}
We shall often simply write $\J_m$.

\begin{proposition}\label{prop:newj}
$\J_m$ is a subalgebra of $\I_m$ and a Banach algebra with the norm
\begin{equation}\label{norm-of-newj}
\|k\|_{\J_m} := \| k \|_m + \|gk\|_{C^*} + \|kg\|_{C^*}\,.
\end{equation}
Moreover $\J_m$ is holomorphically closed in $\I_m$ (and, therefore, in $C^* (X,\F)$).
\end{proposition}

\subsection{Schatten extensions}
Let $(Y,\F)$, $Y:=\tN\times_\Gamma T$,  be a foliated $T$-bundle {\it without} boundary; for example $Y=\pa X\equiv \pa X_0$.
Consider $(\cyl (Y), \F_{\cyl})$ the associated foliated cylinder.
Recall the function
$\chi^0_{{\rm cyl}}$ (often just $\chi^0$), the   function on the cylinder induced by the characteristic function of $(-\infty,0]$ in $\RR$.
Notice that the definition of Schatten norm also apply to $(\cyl (Y), \F_{\cyl})$, viewed as a foliated $T$-bundle with cylindrical ends.
Let $\Psi^{-p}_{\RR,c}(G_{\cyl})\equiv \Psi^{-p}_{c}(G_{\cyl}/\RR_{\Delta})$
be the space of $\RR\times\Gamma$-equivariant families of pseudodifferential operators
of order $-p$ on the fibration $(\RR\times \tN)\times T\to T$ with $\RR\times\Gamma$-compact support.
Consider an element  $\ell\in \Psi^{-p}_{c}(G_{\cyl}/\RR_{\Delta})$; then we know that
$\ell$ defines a bounded operator from the Sobolev field $\E^{k}$ to the Sobolev field $\E^{k+p}$.
See \cite{MN}, Section 3.
Let us denote, as in \cite{MN}, the operator norm of a bounded operator $L$ from $\E^k$ to $\E^j$ as
$\| L \|_{j,k}$; notice the reverse order. For a $\RR\times \Gamma$-invariant, $\RR\times\Gamma$-compactly supported pseudodifferential operator of order $(-p)$, $P$, we consider the norm
\begin{equation}\label{triple-norm}
||| P |||_p \,:= \,{\rm max}(\| P \|_{-n,-n-p}\,,\, \|P\|_{n+p,n})
\end{equation}
with $n$ a fixed integer strictly greater than $\dim N$.
We denote the closure of
$\Psi^{-p}_{c}(G_{\cyl}/\RR_{\Delta})$ with respect to the norm $|||\,\cdot\,|||_p$
by $\operatorname{OP}^{-p}  (\cyl (Y),\F_{\cyl})$. We shall often  write $\operatorname{OP}^{-p}$.

\begin{proposition}\label{prop:b-o-sub-of-b-star}
$\operatorname{OP}^{-p} (\cyl (Y),\F_{\cyl})$ is a Banach algebra and a subalgebra of $B^* (\cyl (Y),\F_{\cyl})$
\end{proposition}

Consider now the bounded linear map $\partial_3^{{\rm max}} : B^*\to \End_{\Gamma} \H$ given
by $\partial_3^{{\rm max}} \ell:= [\chi^0,\ell]$. Consider in $B^*$ the Banach subalgebra
$\operatorname{OP}^{-1} (\cyl (Y),\F_{\cyl})$ and consider in
$\End_{\Gamma} \H$ the subalgebra
$\mathcal{J}_m (\cyl (Y),\F_{\cyl})$. Let  $\partial_3$ be the restriction
of $\partial_3^{{\rm max}} $ to $\operatorname{OP}^{-1} (\cyl (Y),\F_{\cyl})$.
Since $\|\cdot\|\leq |||\,\cdot\,|||$
we see that $\partial_3$ is also bounded.
 Let $\D:=\{\ell\in \operatorname{OP}^{-1} (\cyl (Y),\F_{\cyl})\;\;|\;\; \partial_3 (\ell)\in
\mathcal{J}_m (\cyl (Y),\F_{\cyl})\}$.
One can prove that $\partial_3\,|_{\D}$ induces a closed derivation
$\overline{\delta}_3$ with domain $\D$. This is clearly a closed extension
of the derivation $\delta_3$, $\delta_3 (\ell)=[\chi^0,\ell]$, considered
in Subsection \ref{subsection:sigma3}.
\begin{definition}\label{def:Bp}
If
$m\geq 1$
 we define
$\mathcal{D}_{m}(\cyl (Y), \mathcal F_{{\rm cyl}})$
as $\mathrm{Dom} \,\overline{\delta}_3$
 endowed with norm
\begin{equation}\label{Bp-norm}
\|\ell \|_{\mathcal{D}_{m}} := ||| \ell  ||| + \| [\chi_{{\rm cyl}}^0,\ell]  \|_{\J_m} \,.
\end{equation}
We shall often simply write $\mathcal{D}_{m}$ instead of $\mathcal{D}_{m}(\cyl (Y), \mathcal F_{{\rm cyl}})$.
 \end{definition}

\begin{proposition}\label{prop:shatten-is-ideal} Let
$m\geq 1$,
then
$\mathcal{D}_{m}$ is a Banach algebra with respect to \eqref{Bp-norm}
and a subalgebra of $B^*\equiv B^* (\cyl(Y), \mathcal F_{{\rm cyl}})$. Moreover,
$\mathcal{D}_{m}$ is holomorphically closed  in $B^*$.
\end{proposition}
The Banach algebra we have defined is still too large for the purpose of extending
the eta cocyle. We shall first intersect it with another holomophically closed Banach subalgebra
of $B^*$.

Observe that there exists an action of $\RR$ on $\Psi^{-1}_{c}(G_{\cyl}/\RR_{\Delta})\subset
\operatorname{OP}^{-1} (\cyl (Y),\F_{\cyl})\subset B^*$ defined by
\begin{equation}\label{r-action-alpha}
\alpha_t (\ell):= e^{its} \ell e^{-its}\,,
\end{equation}
with $t\in\RR$, $s$ the variable along the cylinder and $\ell\in
\Psi^{-1}_{c}(G_{\cyl}/\RR_{\Delta})$. Note that $\alpha_t (\ell)$ is again $(\RR\times\Gamma)$-equivariant.
It is clear that $||| \alpha_t (\ell) ||| =  ||| \ell |||$; thus, by continuity,
$\{\alpha_t\}_{t\in\RR}$ yields  a well-defined  action, still denoted  $\{\alpha_t\}_{t\in\RR}$,
of $\RR$ on the Banach algebra $\operatorname{OP}^{-1} (\cyl (Y),\F_{\cyl})$. Note that this action is only strongly continuous.
Let $\partial_\alpha: \operatorname{OP}^{-1}\to \operatorname{OP}^{-1}$ be the
unbounded derivation associated to $\{\alpha_t\}_{t\in\RR}$
\begin{equation}\label{pa-alpha}
\pa_\alpha (\ell):= \lim_{t\to 0} \frac{(\alpha_t (\ell)-\ell)}{t}
\end{equation}
By definition $${\rm Dom} (\partial_\alpha)=\{\ell\in  \operatorname{OP}^{-1} \;|\:  \partial_\alpha (\ell)\;\text{exists in }
\operatorname{OP}^{-1}\}.
$$
One can prove that
the derivation $\partial_\alpha$ is  in fact a {\it closed} derivation.

We endow
${\rm Dom} (\partial_\alpha)$ with the graph norm
\begin{equation}\label{norm-b}
|||\ell |||+ ||| \pa_\alpha (\ell) |||\,.
\end{equation}
It is not difficult to see  that
${\rm Dom} (\partial_\alpha)$  is a Banach algebra with respect to \eqref{norm-b}
and, obviously, a subalgebra of $B^*\equiv B^* (\cyl(Y), \mathcal F_{{\rm cyl}})$; moreover
it is holomorphically closed  in $B^*$.

We can now take the intersection of the  Banach subalgebras
$\mathcal{D}_{m}(\cyl (Y), \mathcal F_{{\rm cyl}})$ and $ \mathrm{Dom}(\pa_{\alpha})$:
$$\mathcal{D}_{m,\alpha}(\cyl (Y), \mathcal F_{{\rm cyl}}):= \mathcal{D}_{m}(\cyl (Y), \mathcal F_{{\rm cyl}})\cap
 \mathrm{Dom}(\pa_{\alpha})$$
and we  endow it with the norm
\begin{equation}\label{Bp-norm-bis-0}
\|\ell\|_{m,\alpha}:= ||| \ell  ||| + \| [\chi_{{\rm cyl}}^0,\ell]  \|_{\J_m} +  |||  \partial_\alpha \ell |||\,.
\end{equation}
Being the intersection of two  holomorphically closed dense subalgebras, also the Banach algebra
$\mathcal{D}_{m,\alpha}(\cyl (Y), \mathcal F_{{\rm cyl}})$ enjoys this property.\\
We are finally ready to define the subalgebra we are interested in. Recall the function
$f_{\cyl} (s,y)=\sqrt{1+s^2}$.

\begin{definition}\label{def:Bp-bis}
If
$m\geq 1$
 we define
\begin{equation}\label{new-Bm}
\mathcal{B}_{m}(\cyl (Y), \mathcal F_{{\rm cyl}}):=
\{\ell\in \mathcal{D}_{m,\alpha}(\cyl (Y), \mathcal F_{{\rm cyl}}) \;|\; [f,\ell]\;\;\text{and}\;\;[f,[f,\ell]]\;\;\text{are bounded}\}\,.
\end{equation}
This will be  endowed with norm
\begin{align*}
\|\ell \|_{\mathcal{B}_{m}} := &
\|\ell\|_{m,\alpha} + 2 \|[f,\ell]\|_{B^*} + \| [f,[f,\ell]]\|_{B^*}\\
=&||| \ell  ||| + \| [\chi_{{\rm cyl}}^0,\ell]  \|_{\J_m} +  |||  \partial_\alpha \ell |||
+2 \|[f,\ell]\|_{B^*} + \| [f,[f,\ell]]\|_{B^*}\,.
\end{align*}

One can
prove that $\B_m (\cyl (Y), \mathcal F_{{\rm cyl}})$  is {\it a holomorphically closed dense subalgebra of $B^*$.}
We shall often simply write $\mathcal{B}_{m}$ instead of $\mathcal{B}_{m}(\cyl (Y), \mathcal F_{{\rm cyl}})$.
 \end{definition}

Let us go back to the foliated bundle with cylindrical end $(X,\F)$.
We now define \begin{equation}
\mathcal{A}_{m}(X,\mathcal F):= \{k\in A^*(X,\mathcal F); \pi(k)\in \mathcal{B}_{m}(\cyl (\pa X), \mathcal F_{{\rm cyl}}), t(k)\in \mathcal{J}_{m}(X,\mathcal F)\}
\end{equation}

\begin{lemma}
$\mathcal{A}_{m}(X,\mathcal F)$ is a subalgebra of $A^*(X,\mathcal F)$.
\end{lemma}

Now we observe that, as   vector spaces,
\begin{equation}\label{direct-sum-cal}
\mathcal{A}_m \cong \mathcal{J}_m \,\oplus \, s (\mathcal{B}_{m})\,.
\end{equation}

Granted this result, we endow $\mathcal{A}_m$ with the direct-sum norm:
\begin{equation}\label{direct-sum-norm}
|| k ||_{\mathcal{A}_m} := || t(k) ||_{m} + || \pi (k) ||_{\mathcal{B}_{m}}
\end{equation}
 Obviously $s$ induces a bounded linear map
$\mathcal{B}_{m} \to \mathcal{A}_m$ of Banach spaces.

\begin{proposition}
($\mathcal{A}_m, || \;\; ||_{\mathcal{A}_m})$ is a Banach algebra. Moreover,
 $ \mathcal{J}_m$ is an ideal in $\mathcal{A}_m$ and there is   a short exact sequence of Banach algebras:
\begin{equation}\label{Shatten exact sequence}
0\rightarrow
 \mathcal{J}_m (X,\mathcal F)\rightarrow \mathcal{A}_m (X;\mathcal F)
\xrightarrow{\pi}\mathcal{B}_{m} (\cyl (\pa X), \mathcal F_{{\rm cyl}}) \rightarrow 0 \,.
\end{equation}
Finally,
$t: A^* (X,\F)\to C^*(X,\F)$ restricts to a bounded
section
$t:  \mathcal{A}_m (X,\F)
 \to \mathcal{J}_m (X,\F)
 $

\end{proposition}

\subsection{Derivations}

In order to extend continuously the cyclic cocycles $\tau_{GV}$ and
$(\tau_{GV}^r, \sigma_{GV})$ we need to take into account the modular
automorphism group, thus decreasing further the size of the short exact sequence $0\rightarrow
 \mathcal{J}_m \rightarrow \mathcal{A}_m \xrightarrow{\pi}\mathcal{B}_{m} \rightarrow 0$.
 Consider the two  derivations $\delta_1$ and $\delta_2$ introduced in Subsection \ref{subsect:absolute-taugv}.
Let us consider first $\delta_2$. Recall the $C^*$-algebra
$C^*_\Gamma (\H)\supset C^*(X,\F)$; it is obtained, by definition, by
closing up the subalgebra $C_{\Gamma,c} (\H)\subset \End_{\Gamma} (\H)$ consisting of those
elements that preserve the continuous field $C^\infty_c (\tV\times T, E)$.
We set
${\rm Dom}\,(\delta_2^{{\rm max}})=\{k\in  C_{\Gamma,c} (\H) \,|\, [\phi,k]\in C^*_\Gamma (\H)\}$
and
$$\delta_2^{{\rm max}}: {\rm Dom}\,(\delta_2^{{\rm max}})\to C^*_\Gamma (\H), \;\;\;
\delta_2^{{\rm max}} (k):= [\phi,k]\,.$$
One can prove that $\delta_2^{{\rm max}}$
is closable. Similarly,  with self-explanatory notation, the bimodule derivation
$$\delta_1^{{\rm max}}:
{\rm Dom}\,(\delta_1^{{\rm max}})\to C^*_\Gamma (\H,\H^\prime), \;\;\;
\delta_1^{{\rm max}} (k):= [\dot{\phi},k]\,,$$
with ${\rm Dom}\,(\delta_1^{{\rm max}}):= \{k\in  C_{\Gamma,c} (\H) \,|\, [\dot{\phi},k]\in
C^*_\Gamma (\H,\H^\prime)\}$
 is closable.
Let $\overline{\delta}^{{\rm max}}_j$ be their respective closures;
thus, for example,
 $$\overline{\delta}^{{\rm max}}_2: {\rm Dom}\,\overline{\delta}^{{\rm max}}_2\subset C^*_{\Gamma} (\H)
 \longrightarrow C^*_\Gamma (\H)$$
and similarly for $\delta_1^{{\rm max}}$. Define now
$$\D_2:=\{a\in {\rm Dom}\,\overline{\delta}^{{\rm max}}_2\cap \J_m (X,\F)\;|\;\overline{\delta}^{{\rm max}}_2
a\in \J_m (X,\F)\}$$
and $\overline{\delta}_2: \D_2\to  \J_m (X,\F)$
as the restriction of $\overline{\delta}^{{\rm max}}_2$ to $\D_2$
with values in  $\J_m (X,\F)$. One can show
that $\overline{\delta}_2$ is a closed derivation.
Define similarly $\D_1$ and the closed derivation $\overline{\delta}_1$.
We set
 \begin{equation}\label{gothic-j}
\mathbf{ \mathbf{ \mathfrak{J}_m } } := \mathcal{J}_m \cap
{\rm Dom} (\overline{\delta}_1)
\cap {\rm Dom} (\overline{\delta}_2) \,.
 \end{equation}
with  ${\rm Dom} (\overline{\delta}_1)=\D_1$ and
${\rm Dom} (\overline{\delta}_2)=\D_2$.

Consider  next $\mathcal{B}_m$; we consider  the  derivations
 $\delta_1:=[\dot{\phi}_{\pa},\;\;]$,  $\delta_2:= [\phi_{\pa},\;\;]$
 on the cylinder $\RR\times\pa X_0$;
 we have already encountered these derivations in
 the definition of the eta cocycle $\sigma_{GV}$;
  see more precisely Definition \ref{def:sigma3}.
 Consider first $\delta_2$.  Define a  closed derivation $\overline{\partial}_2$
by taking the closure of the closable derivation $ \Psi^{-1}_c (G_{\cyl}/\RR_{\Delta})\xrightarrow{\partial_2} B^*$,
with $\partial_2 (\ell):= [\phi_\partial,\ell]$ and
with $ \Psi^{-1}_c (G_{\cyl}/\RR_{\Delta}) $ endowed with the norm $|||\cdot|||$.
One can prove that
 $\overline{\partial}_2|_{\mathfrak{D}_{2}}$, with
 $$
 \mathfrak{D}_{2}=
\{b\in {\rm Dom}(\overline{\partial}_2)\;\;|\;\; \overline{\partial}_2 (b)\in \B_{m}\}$$
is a closed derivation with values in $\mathcal{B}_m$. We set
$\overline{\delta}_2:=\overline{\partial}_2|_{\mathfrak{D}_2}$;
thus ${\rm Dom}(\overline{\delta}_2)=\mathfrak{D}_2$ and
$\overline{\delta}_2:=\overline{\partial}_2|_{\mathfrak{D}_2}$. A similarly definition of
$\overline{\delta}_1$ and
 ${\rm Dom} (\overline{\delta}_1)$ can be given.
We set  \begin{equation}\label{gothic-b}
\mathbf{ \mathbf{ \mathfrak{B}_m } } := \mathcal{B}_m \cap
{\rm Dom} (\overline{\delta}_1)
\cap {\rm Dom} (\overline{\delta}_2)\equiv \mathcal{B}_m \cap \mathfrak{D}_1\cap \mathfrak{D}_2\,.
 \end{equation}
 We endow  $\mathbf{ \mathbf{ \mathfrak{B}_m } } $ with the norm
 \begin{equation}\label{norm-gothic-b}
 \| \ell \|_{\mathbf{ \mathbf{ \mathfrak{B}_m } } }:= \| \ell \|_{\B_m} +  \| \overline{\delta}_1 \ell \|_{\B_m} +
 \|\overline{\delta}_2 \ell \|_{\B_m}
 \end{equation}
  \begin{proposition}\label{prop:gothic-b-hol-closed}
 $ \mathbf{ \mathbf{ \mathfrak{B}_{m} } }$ is holomorphically closed in $B^* (\cyl(\pa X),\F_{\cyl})$.
 \end{proposition}
Finally, we consider the Banach algebra $\A_m (X,\F)$ which is certainly contained
in $C^*_{\Gamma} (\H)$, given that $A_c (X,\F)$ is contained in $C_{\Gamma,c} (\H)$.
Consider again $\overline{\delta}^{{\rm max}}_j$ and restrict it to a derivation
with values in $\A_m (X,F)$:
$$\overline{\delta}_2: \D_2 \to \A_m (X,F)$$
with $\D_2=\{a\in {\rm Dom}\,\overline{\delta}^{{\rm max}}_2\;|\;\overline{\delta}^{{\rm max}}_2 a\in
\A_m (X,F)\}$ and similarly for $\overline{\delta}_1$.  We obtain in this way closed derivations
$\overline{\delta}_1$ and $\overline{\delta}_2$ with domains $ {\rm Dom}\overline{\delta}_1=\D_1$ and
$ {\rm Dom}\overline{\delta}_2=\D_2$.
We set
 \begin{equation}\label{gothic-a}
\mathbf{ \mathbf{ \mathfrak{A}_m } } := \mathcal{A}_m \cap
{\rm Dom} (\overline{\delta}_1)
\cap {\rm Dom} (\overline{\delta}_2)\cap \pi^{-1} (  \mathbf{ \mathbf{ \mathfrak{B}_{m} } } ) \,.
 \end{equation}

  \begin{lemma}\label{lemma:frak-sequence}
 The map $\pi$ sends $ \mathbf{ \mathbf{ \mathfrak{A}_m } }$ into
 $\mathbf{ \mathbf{ \mathfrak{B}_m } }$;  $ \mathbf{ \mathbf{ \mathfrak{J}_m } }$ is an
 ideal in  $\mathbf{ \mathbf{ \mathfrak{A}_m } }$  and we obtain a short exact sequence
 of Banach algebras
 \begin{equation}\label{frak-sequence}
 0\rightarrow \mathbf{ \mathbf{ \mathfrak{J}_m } } \rightarrow \mathbf{ \mathbf{ \mathfrak{A}_m } }\xrightarrow{\pi} \mathbf{ \mathbf{ \mathfrak{B}_m } }\rightarrow 0
 \end{equation}
 The section $s$ and $t$ restricts to bounded
sections $s: \mathbf{ \mathbf{ \mathfrak{B}_m } }
\to  \mathbf{ \mathbf{ \mathfrak{A}_m } }$ and
$t:  \mathbf{ \mathbf{ \mathfrak{A}_m } }
 \to \mathbf{ \mathbf{ \mathfrak{J}_m } }
 $.
 Finally, $\mathbf{ \mathbf{ \mathfrak{J}_m } }$ is holomorphically closed in $C^* (X,\F)$.
\end{lemma}

\subsection{Isomorphism of K-groups}

Let $0\to J\to A\xrightarrow{\pi} B\to 0$ a short exact sequence of Banach algebras. Recall that $K_0 (J):=
K_0 (J^+,J)\simeq \Ker (K_0 (J^+)\to \ZZ)$
and that $K(A^+,B^+)= K(A,B)$.
For the definition of relative K-groups we refer, for example, to  \cite{hr-book}, \cite{lmpflaum}.
Recall that a relative $K_0$-element
for $ A\xrightarrow{\pi} B$
is represented by a  triple $(P,Q,p_t)$ with $P$ and $Q$ idempotents in $M_{k\times k} (A)$
and $p_t\in M_{k\times k} (B)$ a
path of idempotents connecting $\pi (P)$ to $\pi (Q)$.
The excision  isomorphism
\begin{equation}\label{excision-general}
\alpha_{{\rm ex}}: K_0 (J)\longrightarrow K_0 (A,B)
\end{equation}
is given by
\begin{equation*}
\alpha_{{\rm ex}}([(P,Q)])=[(P,Q,{\bf c})]
\end{equation*}
with ${\bf c}$ denoting the constant path.
Consider also $\mathbf{ \mathbf{ \mathfrak{J}_m } }:=\mathcal{J}_m \cap {\rm Dom} (\overline{\delta}_1) \cap {\rm Dom} (\overline{\delta}_2) $ and recall that this is a smooth subalgebra of $C^* (X,\F)$: using also the
excision isomorphism, we obtain
\begin{equation}\label{is-frak}
K_0 (A^*, B^*) \simeq K_0 (C^* (X,\F))\simeq K_0 (\mathbf{ \mathfrak{J}_m })\simeq
K_0 (\mathbf{ \mathfrak{A}_m},  \mathbf{\mathfrak{B}_{m}})\,.
\end{equation}

\subsection{Extended cocycles}\label{subsect:extended-cyclic}
Recall, from general theory, that
$[\tau_{GV}]\in HC^2 (J_c)$ and $[(\tau_{GV}^r, \sigma_{GV})]\in HC^2 (A_c,B_c)$
  can be paired with elements in $K_0 (J_c)$ and $K_0 (A_c, B_c)$ respectively.
  See the proof of our index formula below for the definition of the relative pairing.
   Introduce now the $S^{p-1}$ operation and
  $$S^{p-1}\tau_{GV}=:\tau_{2n}\quad\text{and}\quad (S^{p-1}\tau_{GV}^r,  \frac{3}{2p+1}S^{p-1} \sigma_{GV})=:
  (\tau_{2p}^r, \sigma_{(2p+1)}).$$ We obtain in this way cyclic cocycles and thus classes
   $[\tau_{2p}]\in HC^{2p} (J_c)$ and  $[(\tau_{2p}^r, \sigma_{(2p+1)})]\in HC^{2p} (A_c,B_c)$.

  \begin{proposition}\label{prop:extended-cocycles}
  Let $2n$  equal to the dimension of the leaves in $X=\tilde{M}\times_\Gamma S^1$.
  Then the absolute cocycle $\tau_{2n}$
  extends to a bounded cyclic cocycle on $\mathbf{ \mathfrak{J}_{2n+1} }$ and
 the  eta cocycle $\sigma_{(2n+1)}$
extends to a  bounded
  cyclic cocycle on $\mathbf{\mathfrak{B}}_{2n+1}$.
  \end{proposition}

   \begin{proposition}\label{prop:extended-cocycles-tris}
  Let ${\rm deg} S^{p-1} \tau_{GV}^r=2p>m(m-1)^2-2=m^3-2m^2+m-2$, with $m=2n+1$ and $2n$ equal to the dimension of the leaves
  in $(X,\F)$.
  Then the regularized Godbillon-Vey cochain $S^{p-1} \tau_{GV}^r$, which is by definition
  $\tau_{2p}^r$, extends to a  bounded
  cyclic cochain on  $\mathbf{\mathfrak{A}}_{m}$.
    \end{proposition}

\smallskip
  \noindent
  Summarizing: fix $m=2n+1$, with $2n$ equal to dimension of the leaves and  set
$$\mathfrak{J} := \mathbf{ \mathfrak{J}_m }\,,\quad  \mathfrak{A}:= \mathbf{ \mathfrak{A}_m} \,,\quad
\mathfrak{B}: = \mathbf{\mathfrak{B}_{m}}$$
Using the above two Propositions we see that there are well defined classes
   \begin{equation}\label{higher-cocycles-extended}
   [\tau_{2p}]\in HC^{2p} (\mathbf{\mathfrak{J}})\quad\text{for}\quad 2p\geq 2n
   \end{equation}
    \begin{equation}\label{higher-cocycles-extended-bis}
   [(\tau_{2p}^r, \sigma_{(2p+1)})]\in
   HC^{2p} (\mathbf{\mathfrak{A}}, \mathbf{\mathfrak{B}})\quad \text{for}\quad 2p>m(m-1)^2-2
   \,.
   \end{equation}

\section{$C^*$-index classes. Excision}\label{sec:index}

\subsection{Dirac operators}
We begin with a closed foliated bundle $(Y,\F)$, with $Y=\tN\times_\Gamma T$.
We are also given a $\Gamma$-equivariant complex vector bundle $\we$ on $\tN\times T$,
or, equivalently, a complex vector bundle on $Y$. We assume that $\we$ has a $\Gamma$-equivariant
vertical Clifford structure. We obtain in this way a $\Gamma$-equivariant family of Dirac operators $(D_\theta)_{\theta\in T}$
that will be simply denoted by $D$.
If $(X_0,\F_0)$,
$X_0=\tM\times_\Gamma T$, is a
foliated bundle with boundary, as in the previous sections, then we shall assume the relevant
geometric structures to be of product-type near the boundary.
If $(X,\F)$ is the associated foliated bundle with cylindrical ends, then we shall extend all the structure
in a constant way along the cylindrical ends. We shall eventually assume $\tM$ to be of even dimension, the
bundle $\we$ to be $\ZZ_2$-graded and the Dirac operator to be odd and formally self-adjoint.
We denote by $D^\pa\equiv (D^\pa_\theta)_{\theta\in T}$ the boundary family defined by $D^+$.
This is a $\Gamma$-equivariant family of formally self-adjoint first order elliptic differential operators on a
complete manifold.
We denote by $D^{{\rm cyl}}$ the operator induced by  $D^\pa\equiv (D^\pa_\theta)_{\theta\in T}$
on the cylindrical foliated manifold $(\cyl (\pa X),\F_{{\rm cyl}})$; $D^{{\rm cyl}}$ is  $\RR\times\Gamma$-equivariant.
We refer to \cite{MN} \cite{LPETALE} for precise definitions.
 In all of this section
we shall make the following fundamental

\medskip
\noindent
{\bf Assumption.} There exists $\epsilon>0$ such that $\forall \theta \in T$
\begin{equation}\label{assumption}
L^2-{\rm spec} (D^\pa_\theta) \cap (-\epsilon,\epsilon)= \emptyset
\end{equation}

\medskip
\noindent
For specific examples where this assumption is satisfied, see \cite{LPETALE}. We shall concentrate on the spin-Dirac
case, but it will be clear how to extend the results to general Dirac-type operators.

\subsection{Index class in the closed case}

Let $(Y,\F)$ be a closed foliated bundle.
First, we need to recall how in the closed case we can define an index class $\Ind (D)\in K_* (C^*(Y,\F))$.
There are in fact three equivalent description of $\Ind (D)$, each one with its own
interesting features:
\begin{itemize}
\item the Connes-Skandalis index class, defined by the Connes-Skandalis projector $P_Q$
associated to a pseudodifferential parametrix $Q$ for $D$; $Q$ can be chosen of $\Gamma$-compact support;
\item the Wassermann index class, defined by the Wassermann projector $W_D$;
\item the index class of the graph projection, defined by
the graph projection $e_D$.
\end{itemize}
It is well known
that the  three classes introduced above are equal in $K_0 (C^* (Y,\F))$.

\subsection{The relative index class $\Ind(D,D^\pa)$}\label{subsec:relative}

Let now $(X,\F)$ be a foliated bundle with cylindrical ends. Let $(\cyl (\pa X),\F_{{\rm cyl}})$
the associated foliated cylinder. Recall $ 0\rightarrow C^*(X,\mathcal F)\rightarrow
A^* (X;\mathcal F)\xrightarrow{\pi} B^* (\cyl (\pa X), \mathcal F_{{\rm cyl}})\rightarrow 0$,
the Wiener-Hopf extension
of the $C^*$-algebra of translation invariant operators  $B^* (\cyl (\pa X),\F_{{\rm cyl}})$; see Subsection \ref{subsec:extension}.
We shall be concerned with the K-theory group $K_* (C^*(X,\mathcal F))$
and with the relative group $K_* (A^* (X;\mathcal F),B^* (\cyl (\pa X), \mathcal F_{{\rm cyl}}))$.
We shall write more briefly $0\rightarrow C^*\rightarrow
A^* \xrightarrow{\pi} B^{*} \rightarrow 0$,
and  $K_* (A^*,B^{*})$.
Recall that a relative $K_0$-cycle
for $A^*\xrightarrow{\pi} B^{*}$
is a triple $(P,Q,p_t)$ with $P$ and $Q$ idempotents in $M_{k\times k} (A^*)$
and $p_t\in M_{k\times k} (B^{*})$ a
path of idempotents connecting $\pi (P)$ to $\pi (Q)$.

\begin{proposition}\label{prop:relative-indeces}
Let $(X,\F)$ be a foliated bundle with cyclindrical ends, as above. Consider the Dirac operator on $X$, $D=(D_\theta)_{\theta\in T}$.
Assume    \eqref{assumption}. Then
the graph projection $e_D$ and the Wassermann projection $W_D$  define two  relative
classes in $K_0 (A^*, B^{*})$. These two classes are equal and fix
the {\it relative index class} $\Ind(D,D^\pa)$.
\end{proposition}

The relative classs of Proposition \ref{prop:relative-indeces}
 are more precisely given by the triples \begin{equation}\label{graph-triple}
(e_D, \begin{pmatrix} 0&0\\0&1 \end{pmatrix},p_t) \;\text{ with } \;p_t:= e_{tD^{\cyl}}\;\text{ and }\;
(W_D, \begin{pmatrix} 0&0\\0&1 \end{pmatrix}, q_t)
\;\text{ with } \;q_t:= W_{tD^{\cyl}}\,,
\end{equation}
with $ t\in [1,+\infty].$
The content of the Proposition is that these two triples do define  elements in $K_0 (A^*, B^{*})$
and that these two elements are equal.

\subsection{The index class $\Ind (D)$}\label{subsec:absolute}

Recall the results in \cite{LPETALE} where it is proved that
there is  a well defined parametrix $Q$ for $D^+$,
$Q  D^+=\Id-S_{+}$, $D^+ Q=\Id - S_{-}$,
with remainders $S_\pm$ in  $\KK (\E)\equiv C^*(X,\F)$.
Consequently, there is a well defined Connes-Skandalis projector $P_Q$.
The construction explained in \cite{LPETALE} is an extension to the foliated case
of the parametrix construction of Melrose, with particular care devoted to the non-compactness
of the leaves.
\begin{definition}
The index class associated to a Dirac operator on $(X,\F)$ satisfying assumption \eqref{assumption}
is the Connes-Skandalis index class associated to the Connes-Skandalis projector $P_Q$.
It is denoted by $\Ind (D)\in K_0 (C^*(X,\F) )$.
\end{definition}

\subsection{Excision for index classes}
The following Proposition plays a fundamental role in our approach to higher APS index theory:
\begin{proposition}
Let $D=(D_\theta)_{\theta\in T}$ be a $\Gamma$-equivariant family of Dirac operators
on a foliated manifold with cylindrical ends $X=\tM\times_\Gamma T$. Assume that $\tM$
is even dimensional.
Assume \eqref{assumption}. Then
\begin{equation}\label{ex-of-rel-is-ab}
\alpha_{{\rm ex}}\,(\,\Ind (D)\,)= \Ind (D,D^\pa)
\end{equation}
\end{proposition}

\section{Index theorems}\label{sect:index-theorems}

\subsection{Notation}
From now on we shall fix the dimension of the leaves
of $(X,\F)$, equal to $2n$, and
set
\begin{equation}\label{algebras-no-subscripts}
\mathbf{\mathfrak{J}}:=\mathbf{\mathfrak{J}_{2n+1}}\,,\quad \mathbf{\mathfrak{A}}:=\mathbf{\mathfrak{A}_{2n+1}}\,\quad
\text{and}\quad \mathbf{\mathfrak{B}}:=\mathbf{\mathfrak{B}_{2n+1}}
\end{equation}
so that the short exact sequence in \eqref{frak-sequence}, for $m=2n+1$, is denoted simply as
\begin{equation}\label{sequence-no-subscripts}
0\to \mathbf{\mathfrak{J}}\to \mathbf{\mathfrak{A}}\to \mathbf{\mathfrak{B}}\to 0
\end{equation}
{\it This is the intermediate subsequence, between $0\to J_c\to A_c \to B_c\to 0$
and $0\to C^* (X,\F)\to A^*(X,\F) \to B^* (\cyl(\pa X),\F_{\cyl}) \to 0$, that we have mentioned
in the  introductory remarks in Subsection \ref{subsection:intr-remarks}}.
\subsection{Smooth index classes}\label{subsect:smooth-index}
In  Sections \ref{subsec:relative} and \ref{subsec:absolute}
we stated the existence of two $C^*$-algebraic index classes: the index class and the
{\it relative} index class.
We have also seen in Subsection \ref{subsect:extended-cyclic} that the absolute and relative cyclic
 cyclic cocycles $\tau_{GV}$ and $(\tau_{GV}^r,\sigma_{GV})$  extend
from $J_c$ and $A_c\xrightarrow{\pi_c} B_c$
to the smooth  subalgebras $\mathbf{ \mathfrak{J} }$ and $\mathbf{ \mathfrak{A}} \xrightarrow{\pi}
 \mathbf{\mathfrak{B}}$. In order to make use of the latter information, we need to {\it smooth-out}
 our index classes. This is the content of the following

 \begin{theorem}\label{theo:smooth-indeces}
\item{1)}The Connes-Skandalis projector defines a smooth  index class $\Ind ^s (D)
\in K_0 ( \mathbf{ \mathfrak{J} })$; moreover, if $\iota_* :  K_0 ( \mathbf{ \mathfrak{J} })
\to K_0 (C^*(X,\F))$ is the isomorphism induced by the inclusion $\iota$, then
$\iota_* (\Ind ^s (D))=\Ind (D)$.
\item{2)}The graph projections on $(X,\F)$ and $(\cyl(\pa X),\F_{\cyl})$ define
a smooth relative index class $\Ind^s  (D,D^\pa)
\in K_0 (\mathbf{ \mathfrak{A}},  \mathbf{\mathfrak{B}})$; moreover, if
$\iota_* :  K_0 (\mathbf{ \mathfrak{A}},  \mathbf{\mathfrak{B}})
\to K_0 (A^*,B^*)$ is the isomorphism induced by the inclusion $\iota$,
then $\iota_* (\Ind ^s (D,D^\pa))=\Ind (D,D^\pa).$
\item{3)}Finally,
if $\alpha^s_{{\rm ex}}: K_0 ( \mathbf{ \mathfrak{J} })\rightarrow
K_0 (\mathbf{ \mathfrak{A}},  \mathbf{\mathfrak{B}})$ is the smooth excision isomorphism, then
\begin{equation}\label{smooth-ex}\alpha^s_{{\rm ex}}(\Ind^s (D))=\Ind^s (D,D^\pa)\quad\text{in}\quad K_0 (\mathbf{ \mathfrak{A}},  \mathbf{\mathfrak{B}})\,.
\end{equation}
\end{theorem}

\subsection{The higher APS index formula for the Godbillon-Vey cocycle}

We can now  state  a APS formula for the Godbillon-Vey cocycle.
Let us summarize our geometric data.
We have a foliated bundle with boundary $(X_0,\F_0)$, $X_0=\tM\times_\Gamma T$
with $T=S^1$
We assume that the dimension of $\tM$ is even and that
all our geometric structures (metrics, connections, etc) are of product type near the boundary.
We also
consider $(X,\F)$, the associated foliation with cylindrical ends.
We are given a $\Gamma$-invariant $\ZZ_2$-graded hermitian bundle $\widehat{E}$ on the trivial
fibration $\tM\times T$, endowed with a $\Gamma$-equivariant vertical
Clifford structure. We have a resulting $\Gamma$-equivariant family of Dirac operators $D=(D_\theta)$.

Fix $m=2n+1$, with $2n$ equal to dimension of the leaves and  set as before
$$\mathfrak{J} := \mathbf{ \mathfrak{J}_m }\,,\quad  \mathfrak{A}:= \mathbf{ \mathfrak{A}_m} \,,\quad
\mathfrak{B}: = \mathbf{\mathfrak{B}_{m}}$$
We know that there are well defined  index classes
\begin{equation*}
\Ind^s (D)\in K_0 ( \mathfrak{J})\,,\quad \Ind^s (D,D^\pa)\in K_0 (\mathfrak{A}, \mathfrak{B})\,,
\end{equation*}
the first given in terms of a  parametrix $Q$ and the second given in term of the graph projection
$e_{D}$.
 Proposition \ref{prop:extended-cocycles}  and Proposition \ref{prop:extended-cocycles-tris} imply the existence of the following additive maps:
 \begin{align}\label{pairings}
& \langle \;\cdot\;, [\tau_{2p}] \rangle\,: K_0 (\mathfrak{J})\rightarrow \CC\,,\quad 2p\geq 2n\\
&\langle \;\cdot\; , [(\tau^r_{2p},\sigma_{(2p+1)})] \rangle\,: K_0 (\mathfrak{A}, \mathfrak{B})
\rightarrow \CC\,,\quad 2p>m(m-1)^2-2\,.
\end{align}

\begin{definition}
Let $(X_0,\F_0)$, $X_0=\tM\times_\Gamma S^1$, as above and assume \eqref{assumption}.
The Godbillon-Vey higher index is the number
\begin{equation}\label{gvindex}
\Ind_{GV} (D) := \langle \Ind(D), [\tau_{2n}] \rangle.
\end{equation}
\end{definition}
Notice that, in fact, $\Ind_{GV} (D) := \langle \Ind^s (D), [\tau_{2p}] \rangle$ for each
$p\geq n$.

\medskip
The following theorem is the  main results of this paper:

\begin{theorem}
Let $X_0=\tM\times_\Gamma S^1$ be a foliated bundle with boundary
and let $D:= (D_\theta)_{\theta\in S^1}$ be a $\Gamma$-equivariant
family of  Dirac operators as above. Assume \eqref{assumption} on the boundary family.
Fix  $2p>m(m-1)^2-2$ with $m=2n+1$ and $2n$ equal to the dimension of the leaves. Then
the following two equalities hold
\begin{equation}\label{main-dim3}
\Ind_{GV} (D)= \langle \Ind^s (D,D^{\pa}), [(\tau_{2p}^r,\sigma_{2p+1}] \rangle =
\int_{X_0} {\rm AS}\wedge\omega_{GV} - \eta_{GV}
\end{equation}
 with
 \begin{equation}\label{main-eta}
 \eta_{GV}:= \frac{(2p+1)}{p!}\int_0^{\infty}\sigma_{(2p+1)} ([\dot{p_t},p_t],p_t,\dots,p_t)dt\,,\quad p_t:= e_{tD^{{\rm cyl}}}\,,
 \end{equation}
 defining the {\em Godbillon-Vey eta invariant} of the boundary family and ${\rm AS}$ denoting the form
 induced on $X_0$ by the ($\Gamma$-invariant) Atiyah-Singer form for the fibration $\tM\times S^1\to S^1$
 and the hermitian bundle $\widehat{E}$.
\end{theorem}
Notice that using the Fourier transformation the Godbillon-Vey eta invariant $\eta_{GV}$ does depend only on
 the boundary
family $D^\pa\equiv (D^\pa_\theta)_{\theta\in S^1}$.

\begin{proof}
For notational convenience we set $\tau_{2p}\equiv \tau_{GV}$,
$\tau_{2p}^r \equiv \tau_{GV}^r$
and $\sigma_{(2p+1)}\equiv \sigma_{GV}$. We also write $\alpha_{{\rm ex}}$
instead of $\alpha_{{\rm ex}}^s$.
The left hand side of formula \eqref{main-dim3} is, by definition,
the pairing $\langle [P_Q,e_1], \tau_{GV} \rangle$ with $P_Q$ the Connes-Skandalis
projection and $e_1:= \begin{pmatrix} 0&0\\0&1 \end{pmatrix}$. Here we have also used
the remark that $\Ind_{GV} (D) := \langle \Ind^s (D), [\tau_{2p}] \rangle$ for each
$p\geq n$.
Recall that $\alpha_{{\rm ex}} ( [P_Q,e_1])$
is by definition $[P_Q,e_1,{\bf c}]$, with ${\bf c}$ the constant path with value $e_1$.
Since the derivative of the constant path is equal to zero
and since $\tau^r_{GV}|_{\mathbf{\mathfrak{J}}}=\tau_{GV}$, using the obvious extension
of \eqref{regularized-via-t},
we obtain at once the
crucial relation
\begin{equation}\label{restriction-of-reg}
\langle \alpha_{{\rm ex}} ([P_Q,e_1]), [(\tau^r_{GV},\sigma_{GV})] \rangle
= \langle  [P_Q,e_1], [\tau_{GV}] \rangle \,.
\end{equation}
Now we use the excision formula, asserting that  $\alpha_{{\rm ex}} ( [P_Q,e_1])$
is equal, {\it as a relative class},
to $[e_{D}, e_1, p_t]$ with $p_t:= e_{t D^{{\rm cyl}}}$. Thus
$$\langle [e_{D}, e_1, p_t], [(\tau^r_{GV},\sigma_{GV})] \rangle= \langle  [P_Q,e_1], [\tau_{GV}] \rangle$$
which is the first equality in \eqref{main-dim3} (in reverse order).
Using also the definition of the relative pairing we can summarize our
results so far as follows:
\begin{align*}
\Ind_{GV} (D)&:= \langle \Ind^s (D), [\tau_{GV}] \rangle\\
&\equiv  \langle  [P_Q,e_1], [\tau_{GV}] \rangle\\
&= \langle \alpha_{{\rm ex}} ([P_Q,e_1]), [(\tau^r_{GV},\sigma_{GV})] \rangle \\
&= \langle [e_{D}, e_1, p_t ] , [(\tau^r_{GV},\sigma_{GV})] \rangle \\
&: = \frac{1}{p!} \tau^r_{GV} (e_{D}-e_1)+\frac{(2p+1)}{p!} \int_1^{+\infty}\sigma_{GV} ([\dot{p_t},p_t],p_t,\dots,p_t)dt\\
& \equiv  \frac{1}{p!} \tau^r_{GV} (\widehat{e}_{D})+\frac{(2p+1)}{p!} \int_1^{+\infty}\sigma_{GV} ([\dot{p_t},p_t],p_t,\dots,p_t)dt
\end{align*}
with $\widehat{e}_{D}=(D+\mathfrak{s})^{-1}$.
Notice that the convergence at infinity of the integral 
$\int_1^{+\infty}\sigma_{GV} ([\dot{p_t},p_t],p_t,\dots,p_t)dt$ follows from the fact that the pairing
is well defined.
Replace $D$ by $uD$, $u>0$.
We obtain, after a simple change of variable in the integral,
$$\frac{(2p+1)}{p!} \int_u^{+\infty}\sigma_{GV} ([\dot{p_t},p_t],p_t,\dots, p_t,p_t)dt=-\langle \Ind^s (uD), [\tau_{GV}] \rangle
+\frac{1}{p!}\tau^r_{GV} (\widehat{e}_{uD})$$
But the absolute pairing  $\langle \Ind^s (uD), [\tau_{GV}] \rangle$ in independent of $u$ and of course equal to $\Ind_{GV} (D)$; thus
$$\frac{(2p+1)}{p!} \int_u^{+\infty}\sigma_{GV} ([\dot{p_t},p_t],p_t,\dots,p_t,p_t)dt=-\Ind_{GV} (D)+
\frac{1}{p!}\tau^r_{GV} (\widehat{e}_{uD})$$
 The second summand
of the right hand side can be proved to converge
as $u\downarrow 0$ to $\int_{X_0} {\rm AS}\wedge \omega_{GV}$ (this employs Getzler rescaling
exactly as
in \cite{MN}).
Thus the limit
$$\frac{(2p+1)}{p!} \lim_{u\downarrow 0} \,\int_s^{+\infty}\sigma_{GV} ([\dot{p_t},p_t],p_t,\dots,p_t,p_t)dt$$
exists \footnote{
the situation here is similar to the one for the eta invariant
in the seminal paper of Atiyah-Patodi-Singer; the regularity there is a consequence of  their index theorem}
and  is equal to $\int_{X_0} {\rm AS}\wedge \omega_{GV}-\Ind_{GV} (D)$.
The theorem is proved
\end{proof}

\medskip
\noindent
{\bf Remark.}
The classic Atiyah-Patodi-Singer index theorem in obtained proceeding as above, but pairing
the  index class with the absolute 0-cocycle $\tau_0$ and the relative index class
with the relative  0-cocycle $(\tau_0^r,\sigma_1)$.
Equating the absolute
and the relative pairing, as above, we obtain an index theorem. It can be proved that this is precisely
the APS index theorem on manifolds with cylindrical ends; in other words, the eta-term
we obtain is precisely the APS eta invariant for the boundary operator. The classic APS index theorem
from the point of view of relating pairing was announced by the first author in \cite{AMS}.
This approach is also a Corollary of the main result of the recent preprint of Lesch, Moscovici and Pflaum
\cite{LMP2}, that is,  the computation of the Connes-Chern character of the relative homology
cycle  associated to a Dirac operator on a manifold with boundary in terms of local data
and higher eta cochain for the {\em commutative}
algebra of smooth functions on the boundary (see also \cite{Getzler-b} and \cite{WU}).
Needless to say, the results in \cite{LMP2}
go well beyond the computation of the index; however, they don't have much in common with the
non-commutative results presented in this paper.

\subsection{Eta cocycles}\label{subsection:eta-cocycles}

The ideas explained in the previous sections can be extended to
 general cocycles $\tau_k\in HC^k (C^\infty_c (G, (s^*E)^*\otimes r^*E))$; we simply need
 to require that these cocycles  are in the image of a suitable
 Alexander-Spanier homomorphism since we can then replace integrals with
 regularized integrals in the passage from absolute to relative cocycles.
 This general theory will be treated elsewhere.
Here we only want to comment on the particular case
 of Galois coverings, since this case illustrates very well the general framework.
 In this important example the  techniques of this paper can be used in order to give
 an alternative approach to the higher index theory developed in \cite{LPMEMOIRS}, much more
 in line with the original treatment given by Connes and Moscovici in their fundamental
 paper  \cite{CM}.

We now give a very short treatment of
this important example, assuming a certain familiarity with the seminal work of Connes and Moscovici.
Let  $\Gamma\to \tM\to M$ be a Galois covering with boundary
and let $\Gamma\to \tV\to V$ be the associated covering with cylindrical ends. In the closed case
higher indeces for a $\Gamma$-equivariant Dirac operator on $\tM$ are obtained through Alexander-Spanier
cocycles, so we concentrate directly on these.  Let $\phi$ be an Alexander-Spanier $p$-cocycle; for simplicity
we assume that $\phi$ is the sum of decomposable elements given by the cup product of Alexander-Spanier
1-cochains: $\phi= \sum_i \delta f^{(i)}_1 \cup \delta f^{(i)}_2 \cup \cdots \cup \delta f^{(i)}_p$ where
$f_j^{(i)}:\tM\to \CC$ is continuous.
Here we assume that $\delta f_j^{(i)}$, $\delta f^{(i)}_j (\tm,\tm'):=(f^{(i)}_j (\tm')-
f^{(i)}_j (\tm))$ is  $\Gamma$-invariant with respect to the diagonal action
of $\Gamma$ on $\tM\times \tM$.
This is a non-trivial assumption.
 We
shall omit $\cup$ from the notation.   The cochain $\phi$ is a cocycle (where we recall that for an Alexander-Spanier $p$-cochain given by a continuos
function $\phi: \tM^{p+1}\to \CC$ invariant under the diagonal $\Gamma$-action, one sets $\delta\phi (x_0,x_1,\dots,
x_{p+1}):= \sum_0^{p+1} (-1)^i \phi (x_0,\dots,\hat{x}_i,\dots,x_{p+1})$). Always in the closed case we obtain
a cyclic $p$-cocycle for the convolution algebra $C^\infty_c (\tM\times_\Gamma \tM)$ by setting
\begin{equation}\label{as-cyclic}
\tau_\phi (k_0,\dots,k_p)=\frac{1}{p!} \,  \sum_{\alpha\in \mathfrak{S}_p}\, \sum_i\, {\rm sign} (\alpha) \omega_\Gamma
(k_0 \, \delta^{(i)}_{\alpha (1)} k_1 \cdots \delta^{(i)}_{\alpha (p)} k_p)\,,
\end{equation}
 with $\delta^{(i)}_j k := [k,f^{(i)}_j].$
Notice that $[k,f^{(i)}_j]$ is the $\Gamma$-invariant kernel whose value at
$(\tm,\tm')$ is given by \\$k(\tm,\tm') \delta f^{(i)}_j (\tm,\tm')$ which is by definition $k(\tm,\tm') (f^{(i)}_j (\tm')-
f^{(i)}_j (\tm))$; $\omega_\Gamma$ is as usual given by
$\omega_\Gamma (k)=\int_F \Tr_{\tm} k(\tm,\tm)$, with $F$ a fundamental domain for the $\Gamma$-action.

Pass now to manifolds with boundary and associated manifolds with cylindrical ends. Consider
the small subalgebras $J_c (\tV)$, $A_c (\tV)$, $B_c (\pa \tV\times \RR)$ appearing
in the (small) Wiener-Hopf extension constructed in Subection \ref{subsect:small-dense} (just take $T=$point there).
We write briefly  $J_c$, $A_c $, $B_c$ and $0\to J_c\to A_c \xrightarrow{\pi_c}B_c\to 0$. We adopt
the notation of the previous sections.
Given $\phi$ as above, we can clearly define an {\it absolute} cyclic p-cocycle $\tau_\phi$ on $J_c$.
Next, define the $(p+1)$-linear functional $\psi^r_\phi$ on $A_c$ by replacing the integral in $\omega_\Gamma$ with Melrose' regularized
integral. Consider next the
 cyclic $p$-cochain on $A_c$, call it $\tau_\phi^r (k_0, \dots, k_p)$, defined by 
$$\frac{1}{p+1} \left( \psi^r_\phi(k_0,k_1,\dots,k_p)+ \psi^r_\phi(k_1,\dots,k_p,k_0)+\cdots+
\psi^r_\phi(k_p,k_0,\dots,k_{p-1}) \right)\,.$$
 Finally, introduce the new derivation $\delta^{(i)}_{p+1} (\ell):= [\chi^0,\ell]$
with $\chi^0$ the function on $\pa \tV\times \RR$ induced by the characteristic function
of $(-\infty,0]$. Then
the eta cocycle associated to $\tau_\phi$ is given by
\begin{equation}\label{eta-cm}
\sigma_\phi (\ell_0,\dots,\ell_{p+1})=\frac{1}{(p+1)!}  \sum_{\alpha\in \mathfrak{S}_{p+1}}\, \sum_i\, {\rm sign} (\alpha) \omega_\Gamma
(\ell_0 \, \delta^{(i)}_{\alpha (1)} \ell_1 \cdots \delta^{(i)}_{\alpha (p+1)} \ell_{p+1})
\end{equation}
It should be possible to prove, using the techniques of this paper,  that this is a cyclic $(p+1)$-cocycle for $B_c$ and that $(\tau_\phi^r, \sigma_\phi)$ is a {\it relative}
cyclic $p$-cocycle for the pair $(A_c,B_c)$. $\sigma_\phi$ is, by definition, the {\it eta cocycle} corresponding to
$\tau_\phi$.

Proceeding exactly as above, thus introducing suitable smooth algebras, extending
the cyclic cocycles, smoothing out the index classes and equating the
absolute pairing $\langle \Ind (\tilde{D}), [\tau_\phi] \rangle$ with the relative pairing
$\langle \Ind (\tilde{D},\tilde{D}^\pa), [\tau_\phi^r,\sigma_\phi] \rangle$ one should obtain a higher (Atiyah-Patodi-Singer)-(Connes-Moscovici)
index formula, with boundary correction term given in terms of
$$ \int^\infty_0 \sigma_\phi ( [\dot{p_t},p_t],p_t,\dots,p_t)dt \quad\text{with}\quad p_t:= e_{t \tilde{D}^{\cyl}}$$
A full treatment of the general theory on foliated bundles,
 together with this important particular case will be treated elsewhere.
{\small \bibliographystyle{plain}
\bibliography{proceed-henri-arxiv}
}

\end{document}